\documentclass[a4paper, 10pt]{amsart}

\usepackage[margin=0.75in,headheight=15pt,footskip=25pt]{geometry}

\usepackage{amsmath}
\usepackage{amssymb}
\usepackage{amsthm}
\usepackage{subcaption}

\usepackage{tikz}
\usepackage{pgfplots}
\usetikzlibrary{calc,patterns,angles,quotes}
\usepackage{enumerate}

\usepackage{hyperref}

\hypersetup{colorlinks=true,allcolors=blue}
\usepackage[usenames,dvipsnames]{pstricks}
\usepackage{siunitx,booktabs}
\newtheorem{theorem}{Theorem}[section]
\newtheorem{corollary}{Corollary}[theorem]

\newtheorem{lemma}[theorem]{Lemma}
\newtheorem{definition}{Definition}[section]
\newtheorem{remark}{Remark}[section]
\newtheorem{assumption}{Assumption}[section]

\numberwithin{figure}{section}
\numberwithin{equation}{section}
\newtheorem{algorithm}{Algorithm}

\newcommand{\tribar}{\vert\kern-0.25ex\vert\kern-0.25ex\vert}

\newcommand{\M}{\bold{M}}

\newcommand{\Q}{\bold{T}}
\newcommand{\D}{\bold{D}}

\newcommand{\CC}{\mathcal{C}}

\newcommand{\ddiv}{\normalfont{\text{div}}\,}

\newcommand{\Th}{\mathcal{T}_h}
\newcommand{\Sh}{\mathcal{S}_h}

\newcommand{\N}{N}
\newcommand{\Nh}{\mathcal{N}_h}
\newcommand{\Nhb}{\mathcal{N}_h^{\partial}}
\newcommand{\Ehh}{\mathcal{E}_h}
\newcommand{\Ehb}{\mathcal{E}_h^{\partial}}
\newcommand{\NT}{N_0}
\newcommand{\Zh}{\mathcal{Z}_h}

\newcommand{\Rn}{\mathbb{R}^{\N}}
\newcommand{\R}{\mathbb{R}}
\newcommand{\al}{\boldsymbol{\alpha}}

\newcommand{\be}{\boldsymbol{\beta}}

\newcommand{\rr}{\boldsymbol{r}}

\newcommand{\bff}{\bold{f}}

\newcommand{\x}{\boldsymbol{x}}

\newcommand{\q}{\boldsymbol{q}}
\newcommand{\varthet}{\boldsymbol{\vartheta}}
\newcommand{\brf}{\boldsymbol{\overline{\mathsf{r}}}}
\usepackage{stmaryrd}
\newcommand{\A}{\boldsymbol{A}}

\begin{document}

\title[Error analysis for an AFC scheme for a nonlinear Scalar Conservation Law Using SSP-RK2]{Error analysis of an Algebraic Flux Correction Scheme for a nonlinear Scalar Conservation Law Using SSP-RK2}

\author{Christos Pervolianakis}
\address{Institut für Mathematik, Friedrich-Schiller-Universität Jena, 07743, Jena, Germany}
\email{christos.pervolianakis@uni-jena.de}

\subjclass[2020]{Primary 65M60, 65M15}

\date{\today}

\keywords{finite element method, error analysis, scalar conservation law, inviscid Burger's equation, algebraic flux correction}

\begin{abstract}
We consider  a scalar conservation law with linear and nonlinear flux function on a bounded domain $\Omega\subset{\R}^2$ with Lipschitz boundary $\partial\Omega.$ We discretize the spatial variable with the standard finite element method where we use a local extremum diminishing flux limiter which is linearity preserving. For temporal discretization, we use the second order explicit strong stability preserving Runge--Kutta method. It is known that the resulting fully-discrete scheme satisfies the discrete maximum principle. Under the sufficiently regularity of the weak solution and the CFL condition $k = \mathcal{O}(h^2)$, we derive error estimates in $L^{2}-$ norm for the algebraic flux correction scheme in space and in $\ell^\infty$ in time. We also present numerical experiments that validate that the fully-discrete scheme satisfies the temporal order of convergence of the fully-discrete scheme that we proved in the theoretical analysis.
\end{abstract}

\maketitle
%%-----------------------------
%%      your text
%%-----------------------------

\maketitle

\section{Introduction}
We shall consider the following scalar conservation law where we seek function $u=u(\x,t)$ for $(\x,t)\in{\Omega}\times [0,T],$ satisfying 
\begin{equation}\label{conv_law}
\begin{cases}
u_t  + \ddiv (\bff(u)) = 0 , & \text{in }{{\Omega}}\times [0,T],\\
u = 0, & \text{on }\partial {\Omega}\times[0,T],\\
u(\cdot,0)  = u^0, & \text {in }{{\Omega}},
\end{cases}
\end{equation}
where $\Omega\subset{\R}^2$ is a bounded domain with Lipschitz boundary $\partial\Omega$. We assume that the convective fluxes $\bff = (f_1,f_2)^T$ satisfying the following assumption.
\begin{assumption}\label{assumption:fluxes}
Assume the flux function $\bff$ of \eqref{conv_law} can be written in the form $\bff(u) = \be\,u^{\ell+1},\,\ell=0,1,$ with $\be = \be(\x,t),\,\x\in\Omega,\,t\in [0,T],\,T>0,$ where $\be = (\beta_1,\beta_2)^T,$ with $\ddiv\be \in L^{\infty}(\Omega).$  In case where $\ell = 1,$ we also assume that $\ddiv\be = 0.$
\end{assumption}
The choice of flux function $\bff$ for $\ell=1,$ in the Assumption \ref{assumption:fluxes}, includes also the case where $\be = (1/2,1/2)^T.$ For this flux function, the \eqref{conv_law} is known as the the inviscid Burger equation. It is well known that the solution $u$ of \eqref{conv_law} for the latter flux function, may develop discontinuities (shocks) in finite time, even if the initial data is smooth. This shock formation is due to the steepening of wavefronts, which can cause the gradient of the solution to become unbounded, see, e.g., \cite{dafermos2005}. Hence, the solutions of \eqref{conv_law} are sought in a time interval $[0,T_{\max}],$ where $T_{\max}>0,$ is the maximum time where $u$ is sufficient regular.  \par
It is well known, see, e.g., \cite{dafermos2005}, that the solution of \eqref{conv_law} is positivity preserving, i.e., 
\begin{align}\label{positivity}
u^0(\x) \geq 0,\;\;\x\in\overline{\Omega} \;\;\;\Rightarrow \;\;\; u(\x,t) \geq 0,\;\;\x\in\overline{\Omega},\;t> 0.
\end{align}
Moreover, it can be shown that the maximum principle holds, i.e.,
\begin{align}\label{max_principle}
 \min_{\x\in\overline{\Omega}}u^0 \leq u(\x,t) \leq \max_{\x\in\overline{\Omega}}u^0,\;\;\x\in\overline{\Omega},\,t>0.
\end{align}
The maximum principle is an important property in the analysis of the nonlinear scalar conservation law of form \eqref{conv_law}. It is of great importance to construct numerical methods that satisfy the discrete analogue of these two properties of \eqref{conv_law}. \par
There exists a wide variety of numerical methods for approximating the scalar conservation law \eqref{conv_law} with the same or different boundary conditions or a more general flux function, see e.g., \cite{kucera2016,leveque2012} and the references therein. More specifically, in the context of the discontinuous Galerkin methods, see e.g., \cite{cockburn1996,kucera2019,zhang2004,zhang2010} and alongside with maximum principle limiter, see, e.g., \cite{cockburn1989,zhang2010_2}. In the context of the continuous finite element methods, see, e.g., \cite{burman2010} and alongside with a local extremum diminishing flux limiter that enforces the maximum principle at discrete level, see e.g., \cite{guermond2014, guermond2016,guermond2017, hajduk2023, kuzmin,kuzmin2002, kuzmin2004, lohmann}.\par
The basis for the methods studied is the variational formulation of the model problem, to find function $u(t)\in H_0^1({\Omega}),$ such that,
\begin{align}
(u'(t), v)  - (\bff(u(t)), \nabla v)  = 0, \;\;\;\forall\, v\in H_0^1({\Omega}).\label{weak_u}
\end{align}
The finite element methods studied are based on triangulations $\Th =\{K\}$ of $\Omega,$ with $h = \max_{\Th}\normalfont{diam}(K).$ We use the finite element spaces
\begin{equation}\label{fem_space}
\Sh : = \left\lbrace \chi \in\mathcal{C}(\overline{{\Omega}})\,:\,\chi|_{K} \in \mathbb{P}_1,\;\forall\;K\in \Th,\,\;\text{such that}\;\,\chi = 0\;\;\text{on}\;\;\partial\Omega\right\rbrace.
\end{equation} 
The semi-discrete approximation of the variational problem \eqref{weak_u}, may be written as follows: Find $u_h(t)\in \Sh$, with $u_h(0) = u^0_h\in \Sh$, such that
\begin{align}
(u_{h,t},\chi) - (\bff(u_h) , \nabla \chi) & = 0,\;\;\forall\chi\in \Sh,
\text{ with } u_h(0)  = u_h^0.\label{semi_fem_u_2D}
\end{align}
where $u_h^0\in \Sh$. We may write \eqref{semi_fem_u_2D} in matrix formulation. Let $\Zh = \lbrace Z_j\rbrace_{j=1}^{\N}$ be the set of nodes in the triangulation $\Th$ and $\lbrace \phi_j \rbrace_{j=1}^{\N}\subset \Sh$ the corresponding nodal basis, with $\phi_j(Z_i)=\delta_{ij}.$ Then, we may write $u_h(t)=\sum_{j=1}^{\N}\alpha_j(t)\phi_j$, with $u_h^0=\sum_{j=1}^{\N}\alpha_j^0\phi_j.$ Therefore, the semi-discrete problem \eqref{semi_fem_u_2D} can then be expressed, with $\al = \al(t),$ where
$\al = (\alpha_1, \dots, \alpha_{\N})^T$ as follows,
\begin{align}
\M\al^\prime(t) - \Q_{\al} \al(t) & = \mathbf{0}, \quad\text{ for } t\in[0,T], \text{ with }\al(0)=\al^0,\label{matrix_u}
\end{align}
where $\mathbf{0}$ the zero vector and the matrix $\M=(m_{ij})$ with elements  $m_{ij}=(\phi_i,\phi_j)$ is the usual mass matrix. To define the elements of the matrix due to the flux function $\bff(u_h),$ we need to use a dual notation, i.e., we will often express its coefficients $\tau_{ij}$, $i,j=1,\dots,\N$, as  functions of an element $\psi\in\Sh$, $\tau_{ij}=\tau_{ij}(\boldsymbol{\psi})=\tau_{ij}(\psi)$,  such that $\psi=\sum_j\psi_j\phi_j\in \Sh$ and $\boldsymbol{\psi}= (\psi_1, \dots , \psi_{\N})^T.$ Thus the elements of $\Q_{\boldsymbol{\psi}}=(\tau_{ij})$,  are defined as follows,
\begin{align}\label{T_def}
\tau_{ij}=\tau_{ij}(\boldsymbol{\psi})=\tau_{ij}(\psi) = (\be\,\phi_j\,\psi^\ell , \nabla \phi_i), \quad \text{for } i,j=1,\dots,\N,\;\text{and}\;\ell=0,1.
\end{align}

The solution of \eqref{semi_fem_u_2D}, should respect the discrete analogue of the maximum principle \eqref{max_principle}, i.e.,
\begin{align}\label{semi_discrete_max_principle}
 \min_{\x\in \overline{\Omega}}u^0_h  \leq u_h \leq  \max_{\x\in \overline{\Omega}}u^0_h  \text{  if and only if } \al(t)\in \mathcal{G},\;\text{with}\;\;\mathcal{G} = [\al^{\min,0},\al^{\max,0}],
\end{align}
where with $\al^{\min,0} := \min_{1\leq i\leq \N}\alpha_i^0,\;\al^{\max,0} := \max_{1\leq i\leq \N}\alpha_i^0,$ i.e., the minimum and maximum of the coefficient vector of $u_h^0.$ \par
Since the nodal basis of $\Sh$ is positive, according to \eqref{semi_discrete_max_principle}, the semi-discrete solution $u_h(t)$ of \eqref{semi_fem_u_2D} is non-negative if and only if the coefficient vector $\al(t)$ is non-negative element-wise. \par

A sufficient condition to ensure the discrete maximum principle \eqref{semi_discrete_max_principle} is that the matrix $\M$ be diagonal with positive diagonal elements and $\Q_{\al}$ has non-positive off-diagonal elements, see, e.g., \cite[Theorem 3.14]{kuzmin}, \cite[Chapter 4]{kuzmin}, \cite{kuzmin2002,kuzmin2004}. Since these conditions are non satisfied by the matrices in \eqref{matrix_u}, as describing in \cite{kuzmin,kuzmin2002,kuzmin2004} and the references therein, we first employ the lumped mass method, which results from replacing the mass matrix $\M$ in \eqref{matrix_u} by a diagonal matrix $\M_L$ with elements $m_i = \sum_{i=1}^{\N}m_{ij}$. Then, we add an artificial diffusion operator $\D_{\al}$ so that the off-diagonal elements of $\Q_{\al} + \D_{\al}$ be non-positive. This technique for \eqref{conv_law} can be found in \cite{kuzmin,kuzmin2002,kuzmin2004} and references therein. The resulting semi-discrete scheme is often called \textit{low-order scheme}, since we introduce an error which may manifest the order of convergence. Indeed, some of the elements of the resulting low-order scheme that the artificial diffusion operator canceled may be harmless to maximum principle, so in order to be as much as possible to the initial semi-discrete scheme and not pollute the order of convergence, we may return some of these, either the whole elements or a portion of these, see also \cite{kuzmin} and in references therein. This procedure is called \textit{algebraic flux correction scheme} or \textit{AFC scheme}. To derive the AFC scheme, we  decompose the error we introduced in the low-order scheme, by adding the artificial diffusion operator, into internodal fluxes, see, e.g., \cite{kuzmin,kuzmin2002,kuzmin2004}. This technique for \eqref{conv_law} can be found in \cite{kuzmin,kuzmin2002,kuzmin2004} and references therein. Then we appropriately restore high accuracy in regions where the solution does not violate the maximum principle.  There exists various algorithms to limit the internodal fluxes. We will consider limiters that satisfy the discrete maximum principle and linearity preservation on arbitrary meshes, as the one proposed in \cite{gabriel4}.\par

High-order temporal accuracy can be achieved by using any high-order Runge--Kutta method. One favorable family of Runge--Kutta methods for the scalar conservation law in the form \eqref{conv_law} is the stong stability preserving (SSP) high-order time discretizations, see, e.g., \cite{gottlieb2001,zhang2010_2}, were developed for the solution of the ODE system that results from the semi-discretization of the hyperbolic partial differential equations  with discontinuous solutions. In that cases, this family of methods, guaranties a desired nonlinear or strong stability property that is already satisfied in the semi-discrete level, e.g., the possible oscillations of the solution. More specifically, the ODE system that results from the semi-discretization using a extremum diminishing limiter is of form $\al'(t) = \A(\al(t))\al(t),\,t\geq 0,$ with $\A(\al(t))\in\mathbb{R}^{\N,\N},$ a linear for $\ell = 0$ and nonlinear for $\ell = 1$ square matrix. The ODE system is also satisfies the discrete maximum principle, see \eqref{semi_discrete_max_principle}. If the explicit Euler preserves this property under a CFL condition, then SSP-RK preserves also this property under the same CFL condition, since the SSP-RK methods are based on explicit Euler in the sense that the intermediate stages are convex combination of the explicit Euler. \par

In this work, our purpose is to analyze the stabilized via algebraic flux correction method semi-discrete scheme of \eqref{semi_fem_u_2D}, see in \cite{kuzmin,kuzmin2002,kuzmin2004} and references therein, using a linearity preserving local extremum diminishing limiter, see, e.g., \cite{gabriel3}. The flux function of the scalar conservation law \eqref{conv_law} is as in Assumption \ref{assumption:fluxes}. The fully-discrete scheme results by using high-order explicit methods that are based in explicit Euler, such as the SSP-RK2.  Our analysis of the stabilized schemes is based on the corresponding one employed in \cite{gabriel4}. Since explicit temporal discretizations will be used, the resulting fully-discrete scheme is linear.

We shall use standard notation for the Lebesgue  and Sobolev spaces, namely we denote  $W^m_p=W^m_p(\Omega)$, $H^m=W^{m}_{2}$, $L^p=L^p(\Omega)$, and  with $\|\cdot\|_{m,p}=\|\cdot\|_{W^m_p}$, $\|\cdot\|_{m}=\|\cdot\|_{H^{m}}$,  $\|\cdot\|_{L^p}=\|\cdot\|_{L^p(\Omega)}$, for $m\in\mathbb{N}$ and $p\in[1,\infty]$, the corresponding norms. \par

The fully-discrete schemes we consider approximate $u^n$ by $U^n \in \Sh$ where $u^n=u(\cdot,t^n)$, $t^n=nk$, $n=0,\dots, \NT$ and $\NT\in\mathbb{N}$, $\NT\ge1$, $k=T/\NT$. Under the sufficient regularity assumptions to the solution of \eqref{conv_law} and choosing time step $k$ and spatial step $h,$ such that $k = \mathcal{O}(h^2),$ we derive error estimates  of the form 
\begin{align*}
\|U^n-u^n\|_{L^{2}} & \le C (k^2 + h^{1/2}),\;\;n\geq 0,
\end{align*}
for the explicit second order strong stability preserving Runge--Kutta (SSP-RK2), see, e.g., \cite{gottlieb2001}.\par

The paper is organized as follows: In Section \ref{section:preliminaries} we introduce notation and recall the semi-discrete  \textit{low-order scheme} and the \textit{AFC scheme} for the discretization of \eqref{conv_law} that can be found in \cite{kuzmin,kuzmin2002,kuzmin2004} and references therein. Further, we recall some auxiliary results for the stabilization terms from \cite{chatzipa3}, that we will employ in the analysis that follows and rewrite  the low-order and AFC scheme, as general semi-discrete scheme. In Section \ref{section:fully_discrete}, we discretize the stabilized semi-discrete scheme via AFC method, in time for the flux function of Assumption \ref{assumption:fluxes}, using the second order explicit strong stability preserving Runge--Kutta (SSP-RK2). For a sufficiently smooth solution of \eqref{conv_law} and $k=O(h^{2})$, we derive error estimates in $L^{2}-$norm.  Finally, in Section \ref{section:numerical_results}, we present numerical experiments, validating our theoretical results.

\bigskip
\section{Preliminaries}\label{section:preliminaries}

\subsection{Mesh assumptions}\label{mesh_assumption}

We consider a  family of regular triangulations $\Th=\{K\}$ of a convex bounded domain $\overline\Omega\subset\R^2$.
We will assume that the family $\Th$ satisfies the following assumption.

\begin{assumption}\label{mesh-assumption}
Let $\Th=\{K\}$ be a  family of regular triangulations  of $\overline\Omega$ such that any edge of any $K$ is either a subset of the boundary $\partial\Omega$ or an edge of another $K \in \Th$, and in addition
\begin{enumerate}
\item  $\Th$ is shape regular, i.e, 
there exists a constant $\gamma>0,$ independent of $K$ and $\Th,$ such that 
\begin{equation}\label{shape_regularity}
\frac{h_K}{\varrho_K} \leq \gamma,\quad \forall K\in\Th,
\end{equation}
where $\varrho_K=\text{diam}(B_K)$, and $B_K$ is the inscribed ball in $K$.
\item The  family of triangulations $\Th$ is quasiuniform, i.e., there exists constant $\varrho>0$ such that
\begin{align}\label{quasi-uniformity}
\frac{\max_{K\in\Th}h_K}{\min_{K\in\Th}h_K} \leq \varrho,\quad\forall K\in\Th.
\end{align}
\end{enumerate}
\end{assumption}

\begin{figure}
\centering
\includegraphics[scale=0.7]{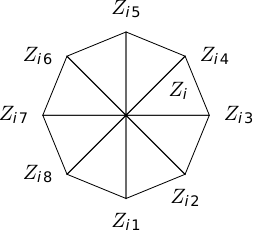} 
\caption{A sub-domain of the triangulation $\Th.$}\label{fig:patches}
\end{figure}

Let $\Nh$ be the the indices of all the nodes of $\Zh$, $\Nh:=\{ i: Z_i \text{ a node of the triangulation } \Th\}$ which  can be splitted  into the indices of the internal nodes, $\Nh^0$,  and  the indices of the nodes  on the boundary $\partial\Omega$,  $\Nhb$, i.e.  $\Nh := \Nh^0\cup \Nhb$.  Also let $\Ehh$ be the set of all edges of the triangulation $ \Th.$ Similarly, we split this set into  the internal edges, $\Ehh^0$  and the edges on the boundary $\partial\Omega$, $\Ehb$, i.e. $\Ehh := \Ehh^{0}\cup\Ehb$. We denote  $\omega_e$ the collection of triangles with a common edge $e\in\Ehh$,  see Fig. \ref{fig:patches},and $\omega_i$, $i\in \Nh$, the collection of triangles with a common vertex $Z_i$, i.e. $\omega_i = \cup_{Z_i\in K}\overline{K},$ see Fig. \ref{fig:patches}. The sets $\Zh(\omega)$ and $\Ehh(\omega)$ contain the vertices and  the edges, respectively, of a 
 subset of $\omega\subset\Th$ and  $\Zh^i$ the set of nodes adjacent to $Z_i$, $\Zh^i:=\{ j: Z_j\in \Zh, \text{adjacent to }Z_i\}$. Using the fact that $\Th$ is shape regular, there exists a constant $C_\gamma$, independent of $h$,  such that the number of vertices in $\Zh^i$  is less than  $C_\gamma$, for $i=1,\dots,\N$. Also $e_{ij}\in\Ehh$ denotes an edge of $\Th$ with endpoints $Z_i$, $Z_j\in \Zh$.

Since $\Th$ satisfies \eqref{quasi-uniformity}, we have for all $\chi\in \Sh,$ cf., e.g., \cite[Chapter 4]{brenner}, 
\begin{equation}\label{eq:inverse_estimate}
\|\chi\|_{L^{\infty}} + \|\nabla \chi\|_{L^{2}}  \le Ch^{-1}\|\chi\|_{L^{2}}.
\end{equation} 

In our analysis, we will employ the following trace inequality which holds for $v\in W^1_p,\;p\in [1,\infty],$ cf. e.g., \cite[Theorem 1.6.6]{brenner}, 
\begin{equation}\label{eq:trace_inequality}
\|v\|_{L^p(\partial\Omega)}  \leq C_\Omega\| v\|_{W^1_p}^{1-1/p}\|v\|_{L^p}^{1/p}.
\end{equation}

\subsection{Stabilized semidiscrete methods}
In order to ensure the discrete maximum principle of the finite element solution, we recall two stabilized schemes, the low order scheme and the algebraic flux correction scheme, see, e.g., \cite{kuzmin,kuzmin2002,kuzmin2004} and references therein.

\subsubsection{Low order scheme}

We replace the composite mass matrix $\M$ by the corresponding lumped mass matrix $\M_L$ and the negative off-diagonal element of $\Q_{\al}$ are cured by an artificial artificial diffusion operator $\D_{\al}=(d_{ij}(\al))$ so that $\Q_{\al}+\D_{\al}\ge0$, element-wise.  To keep our scheme conservative, $\D_{\al}$ must be symmetric with zero row and column sums, cf. \cite{kuzmin,kuzmin2002,kuzmin2004} and references therein.
Also, we will often suppress the index $\al$ in the coefficients $d_{ij}=d_{ij}(\al)$, $i,j=1,\dots,\N$, or express them as  functions of an element  $\psi\in\Sh$, $d_{ij}(\psi)=d_{ij}(\boldsymbol{\psi})$, such that $\psi=\sum_j\psi_j\phi_j\in \Sh$ and  $\boldsymbol{\psi} = (\psi_1, \dots , \psi_{\N})^T$, defined as 
\begin{equation}\label{D_def}
d_{ij} : = \max\{ -\tau_{ij},0, -\tau_{ji}\}=d_{ji}\ge0,\quad\forall j\neq i\ 
\text{ and }\ 
d_{ii} : = -\sum_{j\neq i}d_{ij}.
\end{equation}
For a $\psi\in\Sh,\,\ell=0,1,$ and in view of \eqref{T_def}, we may estimate the elements of the matrix $\D_{\boldsymbol{\psi}}$ as follows,
\begin{equation}\label{est_dij_2D_L_infty}
\begin{aligned}
\vert d_{ij}(\psi)\vert & \leq \vert\tau_{ij}(\psi)\vert + \vert\tau_{ji}(\psi)\vert\\ 
& \leq \|\be\|_{\max}\|\psi^\ell\|_{L^{\infty}}\sum_{K\in \omega_i} (\|\nabla \phi_i\|_{L^{2}(K)}\|\phi_j\|_{L^{2}(K)} + \|\nabla \phi_j\|_{L^{2}(K)}\|\phi_i\|_{L^{2}(K)})\\
& \leq C\|\psi^\ell\|_{L^{\infty}} \sum_{K\in \omega_i}h_K \leq C(\gamma)\|\psi^\ell\|_{L^{\infty}}\,h,
\end{aligned}
\end{equation}
where $\gamma>0$ is the constant of shape regularity, cf. e.g., \eqref{shape_regularity} and for a $\x\in\mathbb{R}^{\N},$ we define $\|\x\|_{\max} = \max_{1\leq i\leq \N}|x_i|.$ 
\begin{remark}\label{remark:q_estimate}
Provided that $\|\psi^\ell\|_{L^{\infty}} \leq M,$ where $M>0$ is uniform and independent of the spatial or temporal discretization, we can conclude that $\vert d_{ij}(\psi)\vert \leq C(\gamma,M)\,h.$ Indeed, while for $\ell = 0$ is always true, we need to ensure it, in the case where $\ell=1.$ 
\end{remark}

The resulting system for the approximation of \eqref{conv_law} is expressed as follows, we seek $\al(t)\in\Rn$ such that, for $t\in[0,T]$,
\begin{align}
\M_L\al^\prime(t) - (\Q_{\al} + \D_{\al}) \al(t) & = 0,\text{ with }\al(0)=\al^0,\label{low_matrix_u}
\end{align}
To write its variational formulation, we define for a function $s\in \Sh,$ its nodal values as $s_i = s(Z_i),\;i=1,\ldots,\N.$ Let for $w\in{\Sh}$, $d_h(w;\cdot,\cdot):{\CC}\times {\CC}\to{\R},$  be a bilinear form defined in \cite{gabriel1}, by
\begin{equation}\label{stab_term}
d_h(w;v,z) := \sum_{i,j=1}^{\N}\,d_{ij}(w)(v_i - v_j)z_i= \sum_{i<j}\,d_{ij}(w)(v_i - v_j)(z_i - z_j),\quad\forall v,z\in{\CC},
\end{equation}
where the last equality is due to the symmetry of matrix $\D,$ see, e.g., \cite{gabriel3}. The bilinear form  $(\cdot,\cdot)_h$ is an inner product in $\Sh$ that approximates $(\cdot,\cdot)$ and is defined by
\begin{equation}\label{quadrature}
(\psi,\chi)_h = \sum_{K\in\Th}Q_h^K(\psi\chi),\ \text{ with }Q_h^K(g) = \frac{1}{3}|K|\sum_{j=1}^3g(Z_j^K)\approx \int_K g\,dx,
\end{equation}
with $\{Z_j^K\}_{j=1}^3$ the vertices of a triangle $K\in\Th.$
In view of \cite{gabriel1}, the algebraic system with the artificial diffusion operator $\D,$ \eqref{low_matrix_u} can be rewritten in the following variational formulation: We seek $u_h(t)\in \Sh$  such that
\begin{align}
(u_{h,t},\chi)_h - (\bff(u_h),\nabla \chi)  + d_h(c_h;u_h,\chi)  & = 0,\ \forall\,\chi\in \Sh,
\label{low_fem_u_2D}
\end{align}
with $u_h(0) = u_h^0.$ The inner product $(\cdot,\cdot)_h$ induces an equivalent norm to $\|\cdot\|_{L^{2}}$ on $\Sh$ where we have the following estimates with constants $C,\,C'$ independent on $h$, such that
\begin{equation}\label{mass_lump_equivalence}
C\|\chi\|_h  \leq \|\chi\|_{L^{2}}  \leq C'\|\chi\|_{h},\ \text{ with }
 \|\chi\|_h = (\chi, \chi)_h^{1/2},\quad\forall\chi\in \Sh.
\end{equation}

\subsubsection{Algebraic flux correction}

The low-order method may harm the convergence rate of the numerical scheme, so following  \cite{kuzmin, kuzmin2004}, one may return some of the canceled fluxed to the semi-discrete scheme \eqref{low_matrix_u} by introducing a flux correction term. Thus, we arrive to the  algebraic flux correction (AFC) scheme, which involves the decomposition of this error into internodal fluxes, which can be used to restore high accuracy in regions where the solution is well resolved and no modifications of the standard FEM are required. There exists various algorithms to implement an AFC scheme. Here we will follow the one proposed in \cite{gabriel1}.\par
Let $\rr = (r_1,\dots,r_\N)^T$  denote  the residual of inserting the operator $\D_{\al}$ in \eqref{matrix_u}, i.e., $\rr(\al) =  \D_{\al}\al.$ Using the zero row sum property of matrix $\D_{\al}$, cf. \eqref{D_def}, we can show, see, e.g., \cite{kuzmin}, that the residual admits a conservative decomposition into internodal fluxes,
\begin{equation}\label{internodal_fluxes_f}
\rr_i = \sum_{j\neq i}\mathsf{r}_{ij},\quad \mathsf{r}_{ji} = -\mathsf{r}_{ij},
\end{equation}
where the amount of mass transported by the raw \textit{antidiffusive flux} $\mathsf{r}_{ij}$ is given by
\begin{align}
\mathsf{r}_{ij} & := \mathsf{r}_{ij}(\al(t)) = d_{ij}(\al(t)) \left(\alpha_i(t) - \alpha_j(t)\right),\;\;\;\;\forall \,j\neq i.\label{flux_D}
\end{align}
The correction terms are defined as
\begin{equation}\label{correction_term}
\overline{\mathsf{r}}_{i}=\sum_{j\neq i}\mathfrak{a}_{ij}\mathsf{r}_{ij},
\end{equation}
where the correction factors $\mathfrak{a}_{ij}=\mathfrak{a}_{ji}\in[0,1],\;i,j=1,\ldots,\N,$ are appropriately defined in view of \eqref{flux_D}.\par
For the rest of this paper we will call the internodal fluxes as anti-diffusive fluxes. Some of these anti-diffusive fluxes are harmless but others may be responsible for the violation of non-negativity. Such fluxes need to be canceled or limited so as to keep the scheme non-negative. Thus, every anti-diffusive flux $\mathsf{r}_{ij}$ is multiplied by a solution-depended correction factor $\mathfrak{a}_{ij}\in[0,1]$, to be defined in the sequel, before it is inserted into the equation. Hence, the AFC scheme is the following: We seek $\al(t)\in\Rn$ such that, for $t\in[0,T]$,
\begin{align}
\M_L\al^\prime(t) - \left(\Q_{\al} + \D_{\al}\right) \al(t) & = \brf\left(\al(t)\right),\,\;\;\text{for}\;t\geq 0\;\text{with}\;\al(0)=\widehat{v},\label{ode_u_afc_2D}
\end{align}
where $\widehat{v},$ is the coefficients vector of $u_h^0\in\Sh.$ \par
To ensure that the AFC scheme maintains satisfies the maximum principle, it is sufficient to choose the correction factors $\mathfrak{a}_{ij}$ such that the sum of anti-diffusive fluxes is constrained by, (cf. e.g., \cite{kuzmin}),
\begin{equation}\label{led_2D}
{Q}_{i}^{-} \leq \sum_{j\neq i}\mathfrak{a}_{ij}\mathsf{r}_{ij} \leq {Q}_{i}^{+},
\end{equation}
and
\begin{equation}\label{alternative_choice_Q}
{Q}_{i}^{+} = q_i(\alpha_i^{\max}(t)-\alpha_i(t))\text{ and }{Q}_{i}^{-} = q_i(\alpha_i^{\min}(t)-\alpha_i(t)),
\end{equation}
and $q_i\geq 0$, $i=1,\dots,\N$, given constants that do not depend on $\al$.
\begin{remark}\label{remark:led}
The criterion \eqref{led_2D} by which the correction factors are chosen, implies that the limiters used in \eqref{ode_u_afc_2D} guarantee that the scheme satisfies the maximum principle. In fact, if $\alpha_i$ is a local maximum, then \eqref{led_2D} implies the cancellation of all positive fluxes. Similarly, all negative fluxes are canceled if $\alpha_i$ is a local minimum. In other words, a local maximum cannot increase and a local minimum cannot decrease. As a consequence, $\mathfrak{a}_{ij}\mathsf{r}_{ij}$ cannot create an undershoot or overshoot at node $i.$ 
\end{remark}
In order to determine the coefficients $\mathfrak{a}_{ij}$, one has to fix first a set of  nonnegative  coefficients $q_{i}$, $i=1,\dots,\N$. In principle the choice  of these parameters $q_{i}$ can be arbitrary. But efficiency and accuracy can dictate a strategy, which does not depend on the fluxes $\mathsf{r}_{ij}$ but on the type of problem ones tries to solve and the mesh parameters. We will not ellaborate more on the choice of $q_{i}$, and for a more detail presentation we refer to \cite{kuzmin} and in the survey \cite{barrenechea2024}. Example of correction factors can be found also in \cite{gabriel1,gabriel4,barrenechea2024,kuzmin,kuzmin2002,kuzmin2004} and the references therein.\par
We shall compute the correction factors $\mathfrak{a}_{ij}$ using Algorithm  \ref{algorithm-1}, which has been proposed   by Kuzmin, cf. \cite[Section 4]{kuzmin} with the choice of $q_i$ as in \cite{gabriel4}.  
\noindent
\begin{algorithm}[Computation of correction factors $\mathfrak{a}_{ij}$] \label{algorithm-1}
Given data: 
\begin{enumerate}
\item The positive  coefficients $q_{i}$, $i,j=1,\dots,\N,$
\item The fluxes $\mathsf{r}_{ij}$, $i\neq j$, $i,j=1,\dots,\N.$
\item The coefficients $\alpha_j,$  $j=1,\dots,\N$.
\end{enumerate}
\noindent
Computation of factors $\mathfrak{a}_{ij},$ for $i,\,j\in\Nh,$ as follows.
\begin{enumerate}
\item
Compute for $i\in\Nh^0,\,j\in\Nh,$ the limited sums $P_{i}^{\pm}:= P_{i}^{\pm}(\al)$ of positive and negative anti-diffusive fluxes
\begin{align*}
P_{i}^{+} = \sum_{j\in \Zh^i}\max\{0, \mathsf{r}_{ij}\},
\quad \text{ and }\quad P_{i}^{-} = \sum_{j\in \Zh^i}\min\{0, \mathsf{r}_{ij}\}.
\end{align*}
\item
Retrieve for $i\in\Nh^0,\,j\in\Nh,$ the local extremum diminishing upper and lower bounds  ${Q}_{i}^{\pm} : = {Q}_{i}^{\pm}(\al),$
\begin{align*}
{Q}_{i}^{+} = q_i(\alpha_i^{\max} - \alpha_i),
\quad\text{ and }\quad{Q}_{i}^{-} = q_i(\alpha_i^{\min} - \alpha_i),
\end{align*}
where $\alpha_i^{\max},\;\alpha_i^{\min}$ are the local maximum and local minimum at $\omega_i.$
\item
Compute for $i\in\Nh^0,\,j\in\Nh,$ also the coefficients $\overline{\mathfrak{a}}_{ij},$ for $j\neq i$ are given by
\begin{equation}\label{correction_factors_definition}
\begin{aligned}
R_{i}^+=\min\left\lbrace 1, \frac{Q_{i}^+}{P_{i}^+} \right\rbrace,\quad R_{i}^-= \min\left\lbrace 1, \frac{Q_{i}^-}{P_{i}^-}\right\rbrace\quad \text{and}\quad \overline{\mathfrak{a}}_{ij} = 
\begin{cases}
R_{i}^+ , &\text{if} \quad\mathsf{r}_{ij} > 0,\\
1, &\text{if} \quad\mathsf{r}_{ij} = 0,\\
R_{i}^-, &\text{if} \quad\mathsf{r}_{ij} < 0.
\end{cases}
\end{aligned}
\end{equation}
\end{enumerate}  
Then, the coefficients $\mathfrak{a}_{ij},$ for $j\neq i$ with $i\in\Nh^0,\,j\in\Nh,$ are given by $\mathfrak{a}_{ij} = \min\{\overline{\mathfrak{a}}_{ij}, \overline{\mathfrak{a}}_{ji}\}$ and $\mathfrak{a}_{ji} = \mathfrak{a}_{ij}.$ For the Dirichlet nodes, i.e., $i\in\Nh^b,\,j\in\Nh^b,$ we set $\mathfrak{a}_{ij} = 1.$
\end{algorithm}
\begin{remark}\label{remark:linearity_preservation1}
In view of \normalfont{\cite[Section 4.2]{gabriel3}}, the above limiter \eqref{correction_factors_definition} is linearity preserving, i.e., for every edge $e$ with endpoints $Z_i,\,Z_j,$
\begin{align}\label{eqn:linear_preserve}
\mathfrak{a}_{ij}(v) = 1,\;\;\text{ if }\;\;v\in \mathbb{P}_1(\mathbb{R}^2),
\end{align}
where $\omega_e$ is the union of triangles with $e$ as a common edges, see the left patch of Fig. \ref{fig:patches}. The linearity preservation is equivalent to have
\begin{align}\label{LP_criterion}
Q_{i}^+ > P_{i}^+\;\;\text{if}\;\;\mathsf{r}_{ij} > 0,\;\;\text{and}\;\;Q_{i}^- < P_{i}^-\;\;\text{if}\;\;\mathsf{r}_{ij} < 0,
\end{align}
see \normalfont{\cite[Section 6]{gabriel4}}.
\end{remark}
\begin{remark}\label{remark:linearity_preservation2}
There exist  $\gamma_i\in\R$, $i\in\Nh^0$,  cf. \normalfont{\cite[Section 6]{gabriel4}}, such that
\begin{equation}\label{def:gamma_i}
v_i - v_i^{\min} \leq \gamma_i (v_i^{\max} - v_i), \quad \forall v\in \mathbb{P}_1(\mathbb{R}^2),
\end{equation}
for $v_i=v(Z_i)$ and $v_i^{\max}$ and $v_i^{\min}$ the local maximum and local minimum, respectively, on $\omega_i$, where $\omega_i$ is the union of triangles with $Z_i$ as a common vertex, see the right patch of Fig. \ref{fig:patches}. Note that in the case of a symmetric $\omega_i$ we have that  $\gamma_i=1$ and $v_i - v_i^{\min} \leq \gamma_i (v_i^{\max} - v_i)$,  for $v\in \mathbb{P}_1(\mathbb{R}^2)$, cf.  \normalfont{\cite[Lemma 6.1]{gabriel4}}.
\end{remark}
\begin{lemma}{\normalfont{\cite[Lemma 6]{gabriel3}}}\label{corrollary:linearity_preservation_limiters}
Let  the positive  coefficients $q_i,$ $i\in\Nh^0$, in Algorithm \ref{algorithm-1} be defined by
\begin{align}\label{def:q_i-1}
q_i : = \gamma_i \sum_{j\in \Zh^i} d_{ij},\quad i\in\Nh^0,
\end{align}
with $\gamma_i$ defined in  \eqref{def:gamma_i}, then the linearity preservation property \eqref{eqn:linear_preserve} is satisfied. 
\end{lemma}

For the correction factors $\mathsf{a}_{ij},\,i,j=1,\ldots,\N,$ that are obtained using Algorithm \ref{algorithm-1} with $q_i,\,i=1,\ldots,\N,$ that satisfies Lemma \ref{corrollary:linearity_preservation_limiters}, we have the following result according to  \cite[Lemma 2.15]{chatzipa3}.
\begin{lemma}{\normalfont{\cite[Lemma 2.15]{chatzipa3}}}\label{lemma:limiters_estimate_localy}
Let the correction factors $\mathfrak{a}_{ij},\,i,j=1,\ldots,\N,$ are obtained using Algorithm \ref{algorithm-1} with $q_i,\,i=1,\ldots,\N,$ that satisfies Lemma \ref{corrollary:linearity_preservation_limiters}. Also for $e_{ij}\in\Ehh^0$, $i\in\Nh^0$, $j\in\Nh$, with endpoints $Z_i, Z_j\in\Zh$, let $\widetilde {\rho}_{ij}(\chi):=\rho_{ij}(\chi)(\chi_i - \chi_j)$,  with $\rho_{ij}(\chi) = \mathfrak{a}_{ij}(\chi)$ or $\rho_{ij}(\chi) = 1 - \mathfrak{a}_{ij}(\chi).$ Then $\widetilde {\rho}_{ij}$ satisfies the following inequality,
\begin{align*}
\vert \widetilde {\rho}_{ij}(\chi) - \widetilde {\rho}_{ij}(\psi)\vert \leq \Lambda_{ij}\sum_{l\in \Zh(\omega_i)}|\chi_l - \psi_l|, \quad \forall \chi,\psi\in \Sh,\;\;i\in\Nh^0,\ j\in\Nh,
\end{align*}
with $\Lambda_{ij} := \Lambda_{ij}(\al,\q) := C(d_{ij}^{-1}(\al)(\max_{1\leq j\leq \N}d_{ij}(\al) + q_i) + 1),$ where the constant $C$ is independent of $h.$ 
\end{lemma}

We can write the algebraic flux correction scheme \eqref{ode_u_afc_2D} in variational formulation. We seek $u_h\in\Sh,$ such that
\begin{equation}\label{semi_fem_u_2D_afc}
\begin{aligned}
(u_{h,t},\chi)_h - (\bff(u_h),\nabla \chi) + \widehat{d}_{h}(u_h;u_h,\chi)  = 0,\;\;\forall\chi\in \Sh,
\text{ with } u_h(0)  = u_h^0,
\end{aligned}
\end{equation}
where the bilinear form $\widehat{d}_{h}(s;\cdot,\cdot)\,:\,{\CC}\times {\CC}\to{\R},$ with  $s\in\Sh$, see \cite{gabriel1}, is defined by, for $v,z\in{\CC}$,
\begin{align}
\widehat{d}_{h}(s;v,z) & := \sum_{i,j=1}^{\N}d_{ij}(s)(1 - \mathfrak{a}_{ij}(s))(v_i - v_j)z_i = \sum_{i<j}d_{ij}(s)(1 - \mathfrak{a}_{ij}(s))(v_i - v_j)(z_i - z_j),\label{stab_term_new2}
\end{align} 
where $v_i := v(Z_i),\;i=1,\ldots,\N,$ for $v\in\CC$ and $\mathfrak{a}_{ij}(s)$ the correction factors that computed in view of Algorithm \ref{algorithm-1} and satisfies Lemma \ref{corrollary:linearity_preservation_limiters}. The last equality holds due to the symmetry of matrix $\D$ and of the coefficients $\mathfrak{a}_{ij},$ see, e.g., \cite{gabriel3}. 

\subsection{Auxiliary results}\label{subsec:auxiliary}
For our analysis, we consider the standard Lagrange interpolation operator $I_h\,:\,\CC(\overline{\Omega}) \to \Sh,$ defined as $I_h\,v(Z_i) = v(Z_i),\,Z_i\in\Zh,$ for a $v\in\CC(\overline{\Omega}).$ The following bounds are hold for every $K\in\Th,$ cf. e.g., \cite[Chapter 4]{brenner} and \cite[Chapter 3]{ciarlet}, 
\begin{align}
\|v - I_hv\|_{L^p(K)} + h_K\|\nabla (v - I_hv)\|_{L^p(K)} &\le Ch^{2}_K\|v\|_{W^{2}_{p}(K)},\quad\forall\,v\in W^{2}_{p},\;\;p\in(1,\infty].\label{I_projection_est_p} 
\end{align}
We also consider the $L^{2}$ projection $P_h:L^{2}\to\Sh$  defined by
\begin{align}
( P_hv - v,& \chi) = 0, \quad\forall \chi\in \Sh.\label{L2_projection_2D}
\end{align}
In view of the mesh Assumption \ref{mesh-assumption}, the projection $P_h$ satisfy the following bounds, cf. e.g., \cite[Chapter 8]{brenner} and \cite{boman2006,crouzeix1987,douglas1975},
\begin{align}
\|P_hv\|_{L^{p}} &\le C\|v\|_{L^{p}},\;\;\;\;\quad\forall v\in L^{p},\;\;p\in[1,\infty]\label{L2_projection_stab}\\
\|P_hv\|_{1,p} &\le C\|v\|_{1,p},\;\;\;\;\quad\forall v\in W^{1}_{p},\;\;p\in[1,\infty]\label{L2_H1_projection_stab}\\
\|v - P_hv\|_{L^{2}} + h\|v - P_hv\|_{1} &\le Ch^2\|v\|_{2},\quad\;\;\,\forall v\in H^2,\label{L2_projection_est2_2D}\\
\|v - P_hv\|_{L^p} &\le Ch^2\|v\|_{2,p},\quad\forall\,v\in W^{2}_{p},\;\;p\in(1,\infty].\label{L2_projection_est_p}
 \end{align}
The inequalities \eqref{L2_projection_stab} and \eqref{L2_H1_projection_stab} can be found in \cite{boman2006,crouzeix1987,douglas1975}. The estimate \eqref{L2_projection_est_p} can be derived using the stability of $L^p-$projection \eqref{L2_projection_stab} together with the interpolation estimates, e.g., \eqref{I_projection_est_p}, for the standard continuous interpolant in $\Sh.$\par
Since the stabilized schemes need a mass lumping, see \eqref{quadrature}, the low-order scheme \eqref{low_fem_u_2D} and the algebraic flux correction scheme (its matrix formulation \eqref{ode_u_afc_2D}, since we have not defined yet its variational formulation), we need to be able to estimate the error of this modification. For the inner product $(\cdot,\cdot)_h$ introduced in \eqref{quadrature}, the following holds.
\begin{lemma}{\normalfont{\cite[Lemma 2.3]{chatzipa2}}}\label{lemma:mass_lump_error}
Let $\varepsilon_{h}(\chi ,\psi): = (\chi , \psi) - (\chi , \psi)_h$. Then,
\begin{equation}
\vert \varepsilon_{h}(\chi ,\psi)\vert\leq Ch^{i+j}\|\nabla^i\chi \|_{L^{2}}\|\nabla^j\psi\|_{L^{2}},\;\;\;\forall\,\chi ,\psi\in \Sh,\;\;\;\text{and}\;\;\;i,j=0,1,\nonumber
\end{equation}
where the constant $C$ is independent of $h.$
\end{lemma}

Next,  we recall various results that will be useful in the analysis that follows. 
Using the following lemma  we have that the bilinear form $d_h$, introduced in \eqref{stab_term}, and hence also $\widehat{d}_h$, defined in \eqref{stab_term_new2}, induces a seminorm on $\CC$.
\begin{lemma}{\normalfont{\cite[Lemma 3.1]{gabriel1}}}\label{lemma:seminorm}
Consider any $\mu_{ij} = \mu_{ji}\geq 0$ for $i,j=1,\dots,\N.$ Then,
\begin{equation}
\sum_{i,j=1}^{\N}v_i\mu_{ij}(v_i - v_j) = \sum_{\substack{i,j=1\\i<j}}^{\N}\mu_{ij} (v_i-v_j)^2\geq 0,\ \forall v_1,\dots,v_{\N}\in\mathbb{R}.\nonumber
\end{equation}
\end{lemma}
Therefore, $d_h(w;\cdot,\cdot):{\CC}\times {\CC}\to{\R},$ with $w\in\Sh$, is a non-negative symmetric bilinear form 
which satisfies the Cauchy-Schwartz's inequality,
\begin{equation} \label{Schwartz_ineq_afc_2D}
|d_h(w;v,z)|^2\leq d_h(w;v,v) d_h(w;z,z),\quad\forall v,z\in{\CC},
\end{equation}  
and thus induces a seminorm on ${\CC}$. 

The bilinear forms $d_h$ and $\widehat d_h$, introduced in \eqref{stab_term} and \eqref{stab_term_new2}, can be written due to symmetry of $\mathfrak{a}_{ij},\,d_{ij},\,i,j=1,\ldots,\N,$ see, e.g., \cite{gabriel3}, as $\overline{d}_h(w,s;\cdot,\cdot):{\CC}\times {\CC}\to{\R},$ with $s\in\Sh$, where
\begin{align}
\overline{d}_{h}(s ;v ,z ) & := \sum_{i<j}\,d_{ij}(s)\rho_{ij}(s)(v_i - v_j)(z_i - z_j),\quad\forall\, v,z\in{\CC},\label{stab_term_afc_2D}
\end{align}
with $\rho_{ij}(s) = \rho_{ji}(s)\in[0,1]$. Note that, for $\rho_{ij}=1$ we have $\overline{d}_{h}=d_h$ and for $\rho_{ij}(s) = 1-\mathfrak{a}_{ij}(s)$, we get $\overline{d}_{h}=\widehat d_{h}.$

For the bilinear form  \eqref{stab_term_afc_2D} the following bound also holds, see \cite[Lemma 2.17]{chatzipa3}. While in \cite{chatzipa3} the proof is based on Neumann boundary conditions, the proof can be easily extended in the case of Dirichlet boundary conditions.
\begin{lemma}{\normalfont{\cite[Lemma 2.17]{chatzipa3}}}\label{lemma:stability_stab_term2}
For $v\in H^2$ and $\rho_{ij}$ satisfies Lemma \ref{lemma:limiters_estimate_localy} there are exists a constant $C$ such that for all $\psi,\,\chi\in \Sh,$
\begin{align}
\vert\overline{d}_{h}(\psi;\psi,\chi)\vert & \leq C(h\|\nabla \psi\|_{L^{\infty}} + \|\q\|_{\max})(\|\nabla (\psi - v)\|^2_{L^{2}} + h^2\|v\|^2_{2})^{1/2}\|\nabla \chi\|_{L^{2}},\label{stability_stab_term_D}
\end{align}
where $\|\q\|_{\max} = \max_{1\leq i\leq \N}|q_i|,$ where $q_i,\;i=1,\ldots,\N,$ are the coefficients on the Algorithm \ref{algorithm-1}, that are used to compute the coefficients $\mathfrak{a}_{ij}(\psi),\,i,j=1,\ldots,\N,$ for the stabilization term \eqref{stab_term_afc_2D}. 
\end{lemma}

\section{Fully-discrete scheme}\label{section:fully_discrete}

For the temporal discretization of \eqref{semi_fem_u_2D_afc}, we will use the second order accurate Strong Stability Preserving Runge--Kutta (SSP-RK), see, e.g., \cite{gottlieb2001}, in uniform partition of the temporal domain. The family of explicit strong stability preserving Runge--Kutta is based on explicit Euler in the sense that the intermediate stages are convex combination of the explicit Euler. Thus, to satisfy the discrete maximum principle, it need to ensure this property for the explicit Euler, see e.g., \cite{burman2015,guermond2017,kuzmin,kuzmin2002,kuzmin2004}. \par
 
Let $\NT\in\mathbb{N}$, $\NT\ge1$, $k=T/\NT$ and $t^n=nk$, $n=0,\dots, \NT$. We seek $U^n\in \Sh$, approximation of $u^n=u(\cdot,t^n)$ for $n=1,\dots,\NT$, such that, 
\begin{equation}\label{fully_disc_1}
\begin{aligned}
(\overline{\partial} U^{n,1},\chi)_h - ( \bff(U^{n-1}) , \nabla \chi) + \widehat{d}_{h}(U^{n-1};U^{n-1}, \chi)  & = 0,
\end{aligned} 
\end{equation}
for $\chi\in\Sh$ and $\overline{\partial} U^{n,1} = (U^{n,1} - U^{n-1})/k$. Moreover, 
\begin{equation}\label{fully_disc_2}
\begin{aligned}
(\widetilde{\partial} U^{n},\chi)_h  - \frac{1}{2}(\bff(U^{n,1}) , \nabla \chi)  + \frac{1}{2}\widehat{d}_{h}(U^{n,1};U^{n,1}, \chi) & = 0,
\end{aligned} 
\end{equation}
for $\chi\in\Sh$ and with $U^0 = u_h^0\in\Sh$ and $\widetilde{\partial} U^n = (U^n - \frac{1}{2} U^{n-1} - \frac{1}{2} U^{n,1})/k$. The resulting fully-discrete scheme \eqref{fully_disc_1}--\eqref{fully_disc_2} is linear. \par
\begin{definition}\label{definition:correction_factors_der}
The correction factors in the stabilization term in \eqref{fully_disc_1}--\eqref{fully_disc_2}, are computed in view of Algorithm \ref{algorithm-1} and satisfies Lemma \ref{corrollary:linearity_preservation_limiters}. More specifically, for a finite element function $\psi\in\Sh,$ with coefficient vector $\varthet\in\mathbb{R}^{\N},$ i.e., $\psi = \sum_{j=1}^{\N}\vartheta_j\phi_j,$ the correction factors $\mathfrak{a}_{ij}(\psi),\,i,j=1,\ldots,\N,$ are computed as follows.
\begin{align*}
\text{The }\mathfrak{a}_{ij}(\psi)\text{ are computed from Algorithm }\ref{algorithm-1}\text{ with }Q^{\pm}(\varthet),\,P^{\pm}(\varthet),\,\text{and }\,q_i = \sum_{j\in\Zh^i}d_{ij}(\psi).
\end{align*}
\end{definition}

The linear system \eqref{fully_disc_1}--\eqref{fully_disc_2} can be also written in matrix formulation. To do this, we introduce the following notation. Let ${\al}^{n,1}=({\alpha}^{n,1}_1,\dots,{\alpha}^{n,1}_{\N})^T,\,{\al}^n=({\alpha}^n_1,\dots,{\alpha}^n_{\N})^T$ the coefficients, with respect to the basis of $\Sh$ of $U^{n,1},\,U^{n}\in\Sh,$ respectively. Then \eqref{fully_disc_1}--\eqref{fully_disc_2} can be written as
\begin{equation}\label{matrix_systemSSP2}
\begin{aligned} 
\M_L\al^{n,1} & = (\M_L + k\,(\Q_{\al^{n-1}} + \D_{\al^{n-1}}  ))\al^{n-1} + k\,\brf(\al^{n-1}) \\
\M_L\al^{n} & = \frac{1}{2}\M_L\al^{n-1} + \frac{1}{2} (\M_L + k\,(\Q_{\al^{n,1}} + \D_{\al^{n,1}}))\al^{n,1} + \frac{k}{2}\,\brf(\al^{n,1}).
\end{aligned}
\end{equation}
It is clear, that the well-posedness of \eqref{matrix_systemSSP2} and as a result of \eqref{fully_disc_1}--\eqref{fully_disc_2} is equivalent to the invertibility of $\M_L,$ which is true for all $k,\,h.$\par
Under suitable smoothness assumptions on the solution of \eqref{conv_law}, we can derive error estimates for the fully-discrete scheme \eqref{fully_disc_1}--\eqref{fully_disc_2}. We will derive error estimates concerning the stabilized fully-discrete schemes. Before we prove the error estimates, let us prove an important a-priori estimate for the finite element solutions $U^{n,1},\,U^n.$

\subsection{Maximum principle}\label{subsection:DMP}

In this section we will discuss known results about the solution $U^n$ of the fully-discrete scheme \eqref{fully_disc_1}--\eqref{fully_disc_2} that satisfies the maximum principle for all $k,\,h.$. Since the family of explicit strong stability preserving Runge--Kutta is based on explicit Euler in the sense that the intermediate stages are convex combination of the explicit Euler, it suffices to prove the maximum principle only for explicit Euler.\par
Since $U^{n,1},\,U^n \in \Sh$ for \eqref{fully_disc_1}--\eqref{fully_disc_2} then, they can be written as a linear combination of the basis functions, i.e., 
\begin{align*}
U^{n,1} = \sum_{i=1}^{\N}\alpha_i^{n,1}\phi_i,\;\;U^{n} = \sum_{i=1}^{\N}\alpha_i^{n}\phi_i,
\end{align*}
for  \eqref{fully_disc_1}--\eqref{fully_disc_2}. The basis functions $\{\phi_i\}_{i=1}^{\N}$ are positive due to construction, therefore, 
\begin{align}\label{discrete_max_principle}
\min_{x\in \overline{\Omega}}u^0_h \leq U^{n,1},\,U^n \leq \max_{x\in \overline{\Omega}}u^0_h \text{  if and only if } \al^{n,1}\in \mathcal{G},\;\text{with}\;\;\mathcal{G} = [\al^{\min,0},\al^{\max,0}],
\end{align}
with $\al^{\min,0} := \min_{1\leq i\leq \N}\alpha_i^0,\;\al^{\max,0} := \max_{1\leq i\leq \N}\alpha_i^0.$ The finite element function $u^0_h$ is an sufficient approximation to $u^0$ onto the finite element space $\Sh$ that preserves the sign of the node values of $u^0.$ Similar for the remaining finite element functions.\par
The proof of the following Theorem, can be found in \cite{kuzmin,kuzmin2002,kuzmin2004} and is based on the criterion \eqref{led_2D}. 
\begin{theorem}\label{theorem:positivity_afc_burger_SSP2} 
Assume the correction factors $\mathfrak{a}_{ij}^n,$ for $i,j=1,\dots,\N,$ computed as in Definition \ref{definition:correction_factors_der}. Then for $\al^{n-1}\in\mathcal{G}$ the coefficients $\al^{n,1},\,\al^n \in \mathcal{G}$ of the solution of \eqref{fully_disc_1}--\eqref{fully_disc_2}.
\end{theorem}

\begin{corollary}\label{corollary:apriori}
Assume the correction factors $\mathfrak{a}_{ij}^n,$ for $i,j=1,\dots,\N,$ computed as in Definition \ref{definition:correction_factors_der}. The following uniform a-priori bounds are hold
\begin{align*}
\|U^{n,1}\|_{L^{\infty}} + \|U^n\|_{L^{\infty}} \leq \|\al^0\|_{\max},\;\;\;n\geq 0.
\end{align*}
\end{corollary}
\begin{proof}
In view of Theorem \ref{theorem:positivity_afc_burger_SSP2}, we obtain
\begin{align*}
\|U^{n-1}\|_{L^{\infty}} & = \left\|\sum_{j=1}^{\N}\alpha_j^{n-1}\phi_j \right\|_{L^{\infty}} \leq \|\al^{n-1}\|_{\max} \left\|\sum_{j=1}^{\N}\phi_j\right\|_{L^{\infty}} \leq \|\al^0\|_{\max}.
\end{align*}

\end{proof}

\subsection{Error estimates}

Under suitable smoothness assumptions on the solution of \eqref{conv_law}, we can derive error estimates for the fully-discrete scheme \eqref{fully_disc_1}--\eqref{fully_disc_2}. We will derive error estimates concerning the stabilized fully-discrete schemes. We follow the ideas of \cite{zhang2004}, by splitting the numerical error into a projection
and a discrete error and then estimate the temporal error. Before the main Theorem, let us assume that appropriate regularity that it will be needed for the derivation of the error estimates.

\begin{assumption}\label{assumption:u_regurarity}
Assume that for the unique solution of scalar conservation law, the following estimates are hold. There exist a constant $M>0,$ such that
\begin{align*}
\|u(t)\|_{2,\infty} + \|u_t(t)\|_2 + \|u_{ttt}(t)\|_{L^{2}} \leq M,\;\;\;0\leq t \leq T_{\max},
\end{align*}
where $T_{\max}>0,$ is the maximum time that the solution of scalar conservation law \eqref{conv_law} with the flux function $\bff$ defined in the Assumption \ref{assumption:fluxes}, has unique solution.
\end{assumption}

Following, we will derive error estimates in $L^2-$norm for the space and $\ell^\infty$ for the time for \eqref{conv_law} following the arguments in \cite[Section 3.1]{burman2010}, \cite{zhang2004}.

\begin{theorem}\label{theorem:error_estimates_2D} 
Let $u$ be the unique, sufficiently smooth solution of \eqref{conv_law}, see Assumption \ref{assumption:u_regurarity}, where the flux function is defined according to the Assumption \ref{assumption:fluxes} and $U^{n,1},\,U^n \in \Sh$ the unique solution of \eqref{fully_disc_1}--\eqref{fully_disc_2} at time level $t=t^n.$ Then, for $k$, $h$ sufficiently small and $k = \mathcal{O}(h^2),$ there exists constant $C>0$, independent of $k, h,$ such that for $n=0,\dots,\NT$, we have
\begin{equation}\label{error_estimate}
\begin{aligned}
\|U^n - u^n\|_{L^{2}} & \leq C(k^2+h^{1/2}  ).
\end{aligned}
\end{equation}
\end{theorem}
\begin{proof}
In view of \cite{burman2010,zhang2004}, we define the function
\begin{align}\label{auxililary_function_for_trunc_error}
w(\x,t) : = u(\x,t) + k\, u_t(\x,t),\;\;\forall\,(\x,t)\in \Omega\times [0,T].
\end{align}
Let the finite element functions $\theta^{n,1},\,\theta^{n}\in \Sh,$ defined $\theta^{n,1}  = U^{n,1} - P_hw^{n-1},\,\theta^{n}  = U^{n} - P_hu^{n}$, and  $\rho^{n,1}  = P_hw^{n-1}-w^{n-1},\,\rho^n  = P_hu^n-u^n$, for $n\ge0$, where $P_h\,:\,L^{2} \to \Sh,$ the usual $L^{2}-$projection defined in \eqref{L2_projection_2D}. \par
The error equation for $\theta^{n,1}$, is
\begin{equation}\label{error_eq_1}
\begin{aligned}
(\overline{\partial}\theta^{n,1}, \chi)_h & + \widehat{d}_{h}(U^{n-1};\theta^{n-1}, \chi) = \mathcal{F}_1(\chi),\;\;\;\;\forall\,\chi\in\Sh,
\end{aligned}
\end{equation}
where the functional $\mathcal{F}_1\,:\,\Sh \to \mathbb{R}$ is defined as
\begin{equation}\label{functional_1}
\begin{aligned}
\mathcal{F}_1(\chi) & = (u_t^{n-1}  - \overline{\partial}P_hw^{n-1}, \chi) + ( \bff(U^{n-1}) - \bff(u^{n-1}) , \nabla \chi)\\
&  - \widehat{d}_{h}(U^{n-1};P_hu^{n-1}, \chi) + \epsilon_h(\overline{\partial}P_hw^{n-1},\chi).
\end{aligned}
\end{equation}
We will estimate its four terms, where for the second term, we will prove distinguish cases according to the definition of the function $\bff,$ see Assumption \ref{assumption:fluxes}. The first term, it can be estimated in view of  \eqref{auxililary_function_for_trunc_error}. We have
\begin{align*}
|(u_t^{n-1}  - \overline{\partial}P_hw^{n-1}, \chi)| \leq Ch^2\sup_{t\in[0,T]}\|u_t(t)\|_{2}\|\chi\|_{L^{2}}  \leq Ch^2\|\chi\|_{L^{2}}.
\end{align*}
For the second term, our aim is to prove the following estimate,
\begin{equation}\label{conv_term_estimates_final}
\begin{aligned}
( \bff(U^{n-1}) -  \bff(u^{n-1}) , \nabla \chi) \leq C(h^2 + \|\theta^{n-1}\|_{L^{2}}^2),
\end{aligned}
\end{equation}
for a constant $C,$ independent of $h$ or $k.$\par
Let as assume first that $\ell = 1$ in the Assumption \ref{assumption:fluxes}, i.e.,  $\bff(u) = \be\,u$ with $\be = \be(\x,t),\,\x\in\Omega,\,t\in [0,T],\,T>0,$ where $\be = (\beta_1,\beta_2)^T,$ with $\ddiv\be \in L^{\infty}(\Omega).$ Then, by using an integration by parts formula and in view of the $\chi = 0$ on $\partial\Omega,$
\begin{equation}\label{conv_term_estimates}
\begin{aligned}
(\bff(U^{n-1}) - \bff(u^{n-1}) , \nabla \chi) & = (\be \, (\theta^{n-1} + \rho^{n-1}), \nabla \chi)\\
& = J_1(\chi) + J_2(\chi).
\end{aligned}
\end{equation}
For $\chi = \theta^{n-1},$ we have
\begin{align*}
J_1(\theta^{n-1}) & = \int_\Omega \be \,\theta^{n-1} \cdot \nabla\theta^{n-1}\,dx = \frac{1}{2}\int_\Omega \be \cdot \nabla (\theta^{n-1})^2\,dx \\
& = \frac{1}{2}\sum_{K\in\Th}\left( -\int_K \ddiv\be  \, (\theta^{n-1})^2\,dx + \int_{\partial K}\be \cdot \nu\,(\theta^{n-1})^2\,ds\right)\\
& = -\frac{1}{2}\int_\Omega \ddiv\be  \, (\theta^{n-1})^2\,dx  \\
& \leq \frac{1}{2}\|\ddiv \be\|_{L^{\infty}}\|(\theta^{n-1})^2\|_{L^1} \leq C\|\theta^{n-1}\|_{L^{2}}^2,
\end{align*}
where we have used an integration by parts formula and the fact that $\theta^{n-1}$ is continuous across the internal edges and zero on the boundary edges, since $\theta^{n-1}\in \Sh.$ Further, using the inverse inequality \eqref{eq:inverse_estimate}, 
\begin{align*}
J_2(\chi) & = \int_\Omega \be\, \rho^{n-1} \cdot \nabla\chi\,dx \leq C\|\rho^{n-1}\|_{L^{2}}\|\nabla \chi\|_{L^{2}} \leq Ch\|\chi\|_{L^{2}}.
\end{align*}
Hence, combining the previous estimates, we can estimate the term \eqref{conv_term_estimates_final} due to convection for $\chi = \theta^{n-1}.$\par

Next, we will prove \eqref{conv_term_estimates_final} for the function $\bff(u) = \be\,u^2,\,u\in \mathbb{R}$ with $\be = \be(\x,t),\,\x\in\Omega,\,t\in[0,T_{\max}],\,T_{\max}>0,$ and $\ddiv\be = 0.$ Notice that this case is with $\ell = 1$ in the Assumption \ref{assumption:fluxes}. The second term of \eqref{functional_1}, can be estimated as follows,
\begin{align*}
J(\chi) := (\be\,(U^{n-1})^2 - \be\,(u^{n-1})^2 , \nabla \chi) & = \int_\Omega \be\, ((U^{n-1})^2 - (u^{n-1})^2) \cdot \nabla \chi\,dx\\
& = \int_\Omega \be\, ( U^{n-1}  -  u^{n-1} )( U^{n-1}  +  u^{n-1}) \cdot \nabla \chi\,dx\\
& = \int_\Omega \be\, (\theta^{n-1} + \rho^{n-1})(U^{n-1}  +  u^{n-1}) \cdot \nabla \chi\,dx.
\end{align*}
Notice that 
\begin{align*}
U^{n-1}  +  u^{n-1} = U^{n-1} - u^{n-1} + 2u^{n-1} = 2u^{n-1} + \theta^{n-1} + \rho^{n-1},
\end{align*}
and then
\begin{align*}
J(\chi) & = (\be\, (\theta^{n-1} + \rho^{n-1})^2 , \nabla \chi) + 2(\be\, u^{n-1}(\theta^{n-1}+\rho^{n-1}), \nabla \chi),\;\;\;\forall\,\chi\in\Sh. 
\end{align*}
By expanding the terms on the right hand side, we have
\begin{align*}
J(\chi) = (\be\,((\theta^{n-1})^2 + 2\rho^{n-1}\theta^{n-1} + (\rho^{n-1})^2), \nabla\chi) + 2(\be\,u^{n-1}(\theta^{n-1} + \rho^{n-1}),\nabla \chi) = J_1(\chi) + J_2(\chi).
\end{align*}
To estimate $J_1,$ we work as follows. We set as in the proof of Theorem \ref{theorem:error_estimates_2D} $\chi = \theta^{n-1},$ and then,
\begin{align*}
J_1(\theta^{n-1}) & = \int_\Omega(\be\,((\theta^{n-1})^2 + 2\rho^{n-1}\theta^{n-1} + (\rho^{n-1})^2)\cdot \nabla\theta^{n-1}\,dx\\
& = \int_\Omega \be\,(\theta^{n-1})^2 \cdot \nabla \theta^{n-1}\,dx + 2\int_\Omega \be\,\rho^{n-1}\theta^{n-1}\cdot \nabla\theta^{n-1} \,dx + \int_\Omega \be\,(\rho^{n-1})^2\cdot \nabla\theta^{n-1}\,dx\\
& = J_1^1 + J_1^2 + J_1^3.
\end{align*}
Using an integration by parts formula as in Theorem \ref{theorem:error_estimates_2D}, for the term \eqref{conv_term_estimates}, we have
\begin{align*}
J_1^1 & =  \int_\Omega \be\,(\theta^{n-1})^2 \cdot \nabla \theta^{n-1}\,dx \\
& = \frac{1}{3}\sum_{K\in\Th}\left( -\int_K \ddiv\be  \, (\theta^{n-1})^3\,dx + \int_{\partial K}\be \cdot \nu\,(\theta^{n-1})^3\,ds\right)\\
& = \int_\Omega \ddiv\be  \, (\theta^{n-1})^3\,dx = 0,
\end{align*}
where we have used an integration by parts formula and the fact that $\theta^{n-1}$ is continuous across the internal edges and zero on the boundary edges, since $\theta^{n-1}\in \Sh.$ The difference with the term in \eqref{conv_term_estimates} is that now the flux vector is divergence free, i.e., $\ddiv \be  = 0,$ see, e.g., Assumption \ref{assumption:fluxes}. \par
Next, using Cauchy-Schwartz inequality and the estimates of the $L^{2}-$projection, cf. e.g., \eqref{L2_projection_est_p}, we get,
\begin{align*}
J_1^2 & = 2\int_\Omega \be\,\rho^{n-1}\theta^{n-1}\cdot \nabla\theta^{n-1} \,dx\\
& \leq 2\|\be\|_{\max} \|\rho^{n-1}\|_{L^{\infty}}\|\theta^{n-1}\|_{L^{2}}\|\nabla\theta^{n-1}\|_{L^{2}}\\
& \leq C\|\rho^{n-1}\|_{L^{\infty}}\|\theta^{n-1}\|_{L^{2}}\|\nabla\theta^{n-1}\|_{L^{2}}\\
& \leq Ch^2\sup_{0\leq t\leq T}\|u(t)\|_{2,\infty}\|\theta^{n-1}\|_{L^{2}}\|\nabla\theta^{n-1}\|_{L^{2}}\\
& \leq Ch\|\theta^{n-1}\|_{L^{2}}^2,
\end{align*}
where $\|\x\|_{\max} = \max_{1\leq i\leq \N}|x_i|,\,\x\in\mathbb{R}^{\N}.$ In the last estimate, we have used the inverse inequality, cf. \eqref{eq:inverse_estimate}. For the last term of $I_1,$ we get by using the \eqref{L2_projection_est_p},
\begin{align*}
J_1^3 & = \int_\Omega \be\,(\rho^{n-1})^2\cdot \nabla\theta^{n-1}\,dx\\
& \leq \|\be\|_{\max} \|\rho^{n-1}\|_{L^{\infty}}\|\rho^{n-1}\|_{L^{2}}\|\nabla\theta^{n-1}\|_{L^{2}}\\
& \leq Ch^4\|\nabla\theta^{n-1}\|_{L^{2}} \leq Ch^3\|\theta^{n-1}\|_{L^{2}} \leq Ch^6 + \|\theta^{n-1}\|_{L^{2}}^2.
\end{align*}
Gathering the estimates for $J_1(\theta^{n-1}),$ we obtain that
\begin{align*}
J_1(\theta^{n-1}) \leq C(h^6 + \|\theta^{n-1}\|_{L^{2}}^2).
\end{align*}
Similar arguments can be used to estimate the term $I_2(\theta^{n-1}).$ More specifically, the latter term can be splitted as follows,
\begin{align*}
J_2(\theta^{n-1}) & = 2\int_\Omega \be\,u^{n-1}(\theta^{n-1} + \rho^{n-1}) \cdot \nabla \theta^{n-1}\,dx\\
& = 2\int_\Omega \be\,u^{n-1}\theta^{n-1} \cdot\nabla \theta^{n-1}\,dx  + 2\int_\Omega\be\, u^{n-1}\rho^{n-1} \cdot \nabla \theta^{n-1}\,dx\\
& = J_2^1 + J_2^2.
\end{align*}  
Similar to the term $J_1^1$ above, we have
\begin{align*}
J_2^1 & = 2\int_\Omega \be\,u^{n-1}\theta^{n-1} \cdot \nabla \theta^{n-1}\,dx \\
& = \sum_{K\in\Th}\left( -\int_K \ddiv(\be\,u^{n-1})  \, (\theta^{n-1})^2\,dx + \int_{\partial K}\be\,u^{n-1} \cdot \nu\,(\theta^{n-1})^2\,ds\right)\\
& = \int_\Omega \ddiv(\be\,u^{n-1})\,(\theta^{n-1})^2\,dx \leq \|\ddiv(\be\,u^{n-1})\|_{L^{\infty}}\|\theta^{n-1}\|_{L^{2}}.
\end{align*}
Next, for the second term,
\begin{align*}
J_2^2 & =  2\int_\Omega\be\, u^{n-1}\rho^{n-1} \cdot \nabla \theta^{n-1}\,dx\\
& \leq \|\be\, u^{n-1}\|_{L^{\infty}}\|\rho^{n-1}\|_{L^{2}}\|\nabla\theta^{n-1}\|_{L^{2}} \leq Ch\|\theta^{n-1}\|_{L^{2}},
\end{align*}
where we have used also the inverse inequality \eqref{eq:inverse_estimate}. Gathering the last two estimates, we get
\begin{align*}
J_2(\theta^{n-1}) \leq C(h^2 + \|\theta^{n-1}\|_{L^{2}}^2).
\end{align*}
All together, we obtain \eqref{conv_term_estimates_final}.\par

For the stabilization term of \eqref{functional_1}, we use \eqref{Schwartz_ineq_afc_2D},
\begin{align*}
- \widehat{d}_{h}(U^{n-1};P_hu^{n-1}, \chi) & \leq \widehat{d}_{h}(U^{n-1};P_hu^{n-1}, P_hu^{n-1})^{1/2}\widehat{d}_{h}(U^{n-1};\chi, \chi)^{1/2}\\
& \leq \frac{1}{2}\widehat{d}_{h}(U^{n-1};P_hu^{n-1}, P_hu^{n-1}) + \frac{1}{2}\widehat{d}_{h}(U^{n-1};\chi, \chi)\\
& \leq Ch\,\|(U^{n-1})^{\ell}\|_{L^{\infty}}\|\nabla P_hu^{n-1}\|_{L^{2}}^2 + \frac{1}{2}\widehat{d}_{h}(U^{n-1};\chi, \chi) \\
& \leq Ch + \frac{1}{2}\widehat{d}_{h}(U^{n-1};\chi, \chi),
\end{align*} 
where in the last estimate we have used that $\|(U^{n-1})^{\ell}\|_{L^{\infty}} \leq \|\be\|_{\max}$ for the case where $\ell=0$ and the a-priori estimate $\|U^{n-1}\|_{L^{\infty}} \leq \|\al^0\|_{\max},$ see Corollary \ref{corollary:apriori} for the case where $\ell=1.$ Notice that later we will set $\chi = \theta^{n-1}$ and the second term on the right hand side will be absorbed with the one in the left hand side of \eqref{error_eq_1}. \par
The last term of the functional represents the error due to mass lumping and can be estimated using Lemma \ref{lemma:mass_lump_error} , using Taylor expansion and the inverse inequality \eqref{eq:inverse_estimate}, i.e.,
\begin{align*}
|\epsilon_h(\overline{\partial}P_hw^{n-1},\chi)|\leq Ch^2\|\overline{\partial}\nabla P_hw^{n-1}\|_{L^{2}}\|\nabla \chi\|_{L^{2}} \leq Ch\|\chi\|_{L^{2}}.
\end{align*}
Thus, in total, in view of the above estimates and the Poincare inequality, we get for $\chi = \theta^{n-1},$
\begin{align}\label{output_theta_n1}
(\overline{\partial}\theta^{n,1},\theta^{n-1})_h + \frac{1}{2}\widehat{d}_{h}(U^{n-1};\theta^{n-1}, \theta^{n-1}) \leq C\left(\|\theta^{n-1}\|_{L^{2}}^2  + h \right).
\end{align}
Notice that $\widehat{d}_h(s,\cdot,\cdot),\,s\in\CC,$ induces a seminorm on $\CC,$ see, e.g., Lemma \ref{lemma:seminorm}, thus the second term in the left hand side of \eqref{output_theta_n1} is non-negative and we can absorb due to its non-negativity.\par 
Now, we need to derive a similar estimate for $\theta^{n}.$ The error equation for $\theta^{n},$ satisfies the following  error equation, 
\begin{equation}\label{error_eq_2}
\begin{aligned}
(\widetilde{\partial}\theta^n, \chi)  & + \frac{1}{2}\widehat{d}_{h}(U^{n,1};\theta^{n,1}, \chi) = \mathcal{F}_2(\chi),\;\;\;\;\forall\,\chi\in\Sh,
\end{aligned}
\end{equation}
with $\chi\in \Sh$ and the functional $\mathcal{F}_2\,:\,\Sh\to \mathbb{R}$ defined by
\begin{equation}\label{functionaL2}
\begin{aligned}
\mathcal{F}_2(\chi) & = \frac{1}{2}( w_t^{n-1} - 2\widetilde{\partial}P_hu^n, \chi) + \frac{1}{2}(\bff(U^{n,1}) - \bff(w^{n-1}) , \nabla\chi) \\
& -  \frac{1}{2}\widehat{d}_{h}(U^{n,1};P_hw^{n-1}, \chi)  + \epsilon_h(\widetilde{\partial}P_hu^{n},\chi).
\end{aligned}
\end{equation}
The terms on the right hand side, may be estimated by the following arguments. By the definition of the function $w,$ see \eqref{auxililary_function_for_trunc_error} and using Taylor expansion, we get
\begin{align*}
k\,w_t^{n-1} - 2u^n + u^{n-1} + w^{n-1} = -2\frac{k^3}{3!}u_{ttt}^{n-1} + \mathcal{O}(k^4),
\end{align*}
thus, using elementary calculations and the definition of the $L^{2}-$projection, see \eqref{L2_projection_2D}, we get
\begin{align*}
|( w_t^{n-1} - 2\widetilde{\partial}P_hu^n, \chi)| \leq C(k^2 + h^2)\|\chi\|_{L^{2}}. 
\end{align*}
Next, for the second term of \eqref{functionaL2}, using arguments as before for the estimation of \eqref{conv_term_estimates_final}, we get for $\chi = \theta^{n,1},$ that
\begin{align*}
\frac{1}{2}(\bff(U^{n,1}) - \bff(w^{n-1}) , \nabla\chi)  \leq C(\|\theta^{n,1}\|_{L^{2}}^2 + h^2).
\end{align*}
Further, using similar arguments as in the previous error equation, the stabilization terms can be estimated as
\begin{align*}
-\frac{1}{2}\widehat{d}_{h}(U^{n,1};P_hw^{n-1}, \chi) \leq Ch + \frac{1}{4}\widehat{d}_{h}(U^{n,1};\chi, \chi).
\end{align*}
The last term of the functional represents the error due to mass lumping and can be estimated using Lemma \ref{lemma:mass_lump_error}, using Taylor expansion and the inverse inequality \eqref{eq:inverse_estimate}, i.e.,
\begin{align*}
|\epsilon_h(\widetilde{\partial}P_hu^{n},\chi)|\leq Ch^2\|\widetilde{\partial}\nabla P_hu^{n}\|_{L^{2}}\|\nabla \chi\|_{L^{2}} \leq Ch\|\chi\|_{L^{2}}.
\end{align*}
Thus, in total, in view of the above estimates, we get for $\chi = \theta^{n,1},$
\begin{align}\label{output_theta_n}
(\widetilde{\partial}\theta^{n},\theta^{n,1})_h + \frac{1}{4}\widehat{d}_{h}(U^{n,1};\theta^{n,1}, \theta^{n,1})  \leq C\left(\|\theta^{n,1}\|_{L^{2}}^2 + k^4 + h  \right).
\end{align}
Combining \eqref{output_theta_n1} with \eqref{output_theta_n}, by adding them and multiply by $k,$ we get
\begin{equation}\label{error_eq_final}
\begin{aligned}
\|\theta^{n}\|_{h}^2 - \|\theta^{n-1}\|_{h}^2 & + \frac{1}{2}\widehat{d}_{h}(U^{n-1};\theta^{n-1}, \theta^{n-1}) + \frac{1}{2}\widehat{d}_{h}(U^{n,1};\theta^{n,1}, \theta^{n,1}) \\
& \leq   \|\theta^n - \theta^{n,1}\|_{h}^2  + Ck\left( k^4 + h  +  \|\theta^{n,1}\|_{L^{2}}^2 + \|\theta^{n-1}\|_{L^{2}}^2 \right).
\end{aligned}
\end{equation}
Thus, we need to derive an estimate for $\|\theta^n - \theta^{n,1}\|_{h}.$ To do this, first notice that
\begin{align*}
2k\widetilde{\partial}\theta^{n} - k\overline{\partial}\theta^{n,1} & = \left( (2\theta^n - \theta^{n-1} - \theta^{n,1}) - (\theta^{n,1} - \theta^{n-1})\right)\\
& = 2k (\theta^n -\theta^{n,1}),
\end{align*}
thus to estimate  $\|\theta^n - \theta^{n,1}\|_{h}^2,$ similar to \cite{burman2010}, we multiply \eqref{error_eq_2} by $2$ and we subtract \eqref{error_eq_1} to get
\begin{equation}\label{error_equation_aux2}
\begin{aligned}
(\theta^n - \theta^{n,1}, \chi)_h & = k\,(w_t^{n-1} - \widetilde{\partial}P_hu^{n}, \chi) -  k\,(u_t^{n-1} - \overline{\partial}P_hw^{n-1}, \chi) \\
& + k\,( ( \bff(U^{n,1}) -  \bff(w^{n-1})) - (\bff(U^{n-1}) - \bff(u^{n-1}))  , \nabla\chi)\\
& - k\,( \widehat{d}_{h}(U^{n,1};U^{n,1},\chi) - \widehat{d}_{h}(U^{n-1};U^{n-1} , \chi))\\
& + k\, \epsilon_h(2\widetilde{\partial}P_hu^{n} - \overline{\partial}P_hw^{n-1},\chi)\\
& = I_1(\chi) + I_2(\chi) + I_3(\chi) + I_4(\chi),
\end{aligned}
\end{equation}
where for the first, in view of $L^{2}$ estimates, see \eqref{L2_projection_est2_2D}, and \eqref{auxililary_function_for_trunc_error} and Taylor expansion, we have
\begin{align*}
I_1(\chi) \leq Ck(k^2 + h^2)\|\chi\|_{L^{2}}.
\end{align*}
For second term, for both choices of $\bff,$ we do not use the estimates for \eqref{conv_term_estimates}, but only the Cauchy-Schwartz inequality and the inverse inequality \eqref{eq:inverse_estimate}. More specifically, for $\bff(u) = \be\,u,$ we have
\begin{align*}
I_2(\chi) & \leq Ck\left(h + h^{-1}(\|\theta^{n-1}\|_{L^{2}} +  \|\theta^{n,1}\|_{L^{2}})\right)\|\chi\|_{L^{2}},
\end{align*}
and also for $\bff(u) = \be\,u^2,$ in view of Corollary \ref{corollary:apriori} and the inverse estimate \eqref{eq:inverse_estimate}, we get
\begin{align*}
I_2(\chi) & \leq C\left(\|\theta^{n,1}+\rho^{n,1}\|_{L^{2}}\|w^{n-1} + U^{n,1}\|_{L^{\infty}} + \|\theta^{n-1}+\rho^{n-1}\|_{L^{2}}\|u^{n-1} + U^{n-1}\|_{L^{\infty}}\right)\|\nabla \chi\|_{L^{2}}\\
& \leq Ch^{-1}\left(\|\theta^{n,1}+\rho^{n,1}\|_{L^{2}} + \|\theta^{n-1}+\rho^{n-1}\|_{L^{2}}\right)\|\chi\|_{L^{2}},
\end{align*}
where the latter constant $C$ depends on $\|\al^0\|_{\max}$ and $\sup_{0\leq t\leq T}\|u(t)\|_{L^{\infty}}.$ Further, using $L^{2}$ estimates, see \eqref{L2_projection_est2_2D}, we obtain
\begin{align*}
I_2(\chi) \leq C_1\left(h + h^{-1}\|\theta^{n,1}\|_{L^{2}} + h^{-1}\|\theta^{n-1}\|_{L^{2}}\right)\|\chi\|_{L^{2}},
\end{align*}
where the latter constant $C_1$ depends on constant $C_1$ and in $\sup_{0\leq t\leq T}\|u(t)\|_{2}.$

In view of Lemma \ref{lemma:stability_stab_term2}, the choice of the correction factors, see Remark \ref{remark:linearity_preservation1}, Lemma \ref{corrollary:linearity_preservation_limiters}, the a-priori estimate in Corollary \ref{corollary:apriori} and the inverse estimate \eqref{eq:inverse_estimate}, the remaining terms, can be estimated as
\begin{align*}
I_3(\chi)  & \leq Ck(h\|U^{n,1}\|_{L^{\infty}} + \|\q^1\|_{\max}) (\|\nabla \theta^{n,1}\|_{L^{2}} + h\|P_hw^{n-1}\|_{2})\|\nabla\chi\|_{L^{2}}\\
&  + Ck(h\|U^{n-1}\|_{L^{\infty}}  + \|\q^2\|_{\max})(\|\nabla \theta^{n-1}\|_{L^{2}} + h\|P_hu^{n-1}\|_{2}) \|\nabla\chi\|_{L^{2}}\\
&\leq Ck\,(h + \|\nabla \theta^{n-1}\|_{L^{2}} + \|\nabla \theta^{n,1}\|_{L^{2}} )\|\chi\|_{L^{2}},
\end{align*}
where $\|\q^\ell\|_{\max} = \max_{1\leq i\leq \N}|q_i^\ell|,\,\ell=1,2,$ with $q_i^1 = \sum_{j\in\Zh}d_{ij}(U^{n,1}),\;q_i^2 = \sum_{j\in\Zh}d_{ij}(U^{n-1})$ and thus in view of the estimate \eqref{est_dij_2D_L_infty} and the fact that $\Th$ is shape regular, we obtain that $ \|\q^\ell\|_{\max} \leq C_M\,h,\;\ell=1,2.$

In view of the inverse inequality and the space and time mesh restriction,
\begin{align*}
I_3(\chi) & \leq Ck(h + h^{-1}(\|\theta^{n-1}\|_{L^{2}} + \|\theta^{n,1}\|_{L^{2}}))\|\chi\|_{L^{2}}\\
& \leq Ck^2h^2 + Ck^2h^{-2}(\|\theta^{n-1}\|_{L^{2}}^2 + \|\theta^{n,1}\|_{L^{2}}^2) + \frac{1}{4}\|\chi\|_{L^{2}}^2\\
& \leq Ck^2h^2 + Ck(\|\theta^{n-1}\|_{L^{2}}^2 + \|\theta^{n,1}\|_{L^{2}}^2) + \frac{1}{4}\|\chi\|_{L^{2}}^2,
\end{align*}
where the last inequality holds for $k = \mathcal{O}(h^2).$ For the error due to the mass lumping, in view of Lemma \ref{lemma:mass_lump_error}, by Taylor expansion and the inverse inequality \eqref{eq:inverse_estimate}, we have
\begin{align*}
I_4(\chi) \leq Ch^2(\|\widetilde{\partial}\nabla P_hu^{n}\|_{L^{2}} + \|\overline{\partial}\nabla P_hw^{n-1}\|_{L^{2}})\|\nabla \chi\|_{L^{2}} \leq Ch\|\chi\|_{L^{2}}.
\end{align*}
Setting $\chi = \theta^{n} - \theta^{n,1},$ into \eqref{error_equation_aux2}, and in view of above estimates, we get
\begin{align*}
\|\theta^n - \theta^{n,1}\|_{h}^2 \leq Ck\left(k^2 + h + \|\theta^{n-1}\|_{L^{2}} + \|\theta^{n,1}\|_{L^{2}}\right).
\end{align*}
Inserting this estimate into \eqref{error_eq_final}, we have
\begin{align}\label{estimate_n_prior_final}
\|\theta^{n}\|_{h}^2  \leq (1 + Ck)\|\theta^{n-1}\|_{h}^2 + Ck(k^4 + h +\|\theta^{n,1}\|_{L^{2}}^2).
\end{align}
To conclude to an estimate for $\theta^{n,1}$ from above inequality, we need to derive a sufficient bound for $\|\theta^{n,1} - \theta^{n-1}\|_{h}.$ Recall the error equation \eqref{error_eq_1},
\begin{equation}\label{error_eq_1_2}
\begin{aligned}
(\overline{\partial}\theta^{n,1}, \chi)_h & = (u_t^{n-1}  - \overline{\partial}P_hw^{n-1}, \chi) + (\bff(U^{n-1}) - \bff(u^{n-1}), \nabla\chi)  \\
& - \widehat{d}_{h}(U^{n-1};U^{n-1}, \chi) + \epsilon_h(\overline{\partial}P_hw^{n-1}, \chi)\\
& = I_1(\chi) + I_2(\chi) + I_3(\chi) + I_4(\chi).
\end{aligned}
\end{equation}
Then, we have using the Cauchy-Schwartz inequality,
\begin{align*}
I_1(\chi) = (u_t^{n-1}  - \overline{\partial}P_hw^{n-1}, \chi) \leq Ch^2\sup_{t\in[0,T]}\|u_t(t)\|_{2}\|\chi\|_{L^{2}}  \leq Ch^2\|\chi\|_{L^{2}}.
\end{align*}
Further, similar to the previous estimation for $I_2$ in \eqref{error_equation_aux2}, we use only the Cauchy-Schwartz inequality, to get
\begin{align*}
I_2(\chi) & \leq Ch^{-1}(\|\theta^{n-1}\|_{L^{2}} + \|\rho^{n-1}\|_{L^{2}})\|\chi\|_{L^{2}}.
\end{align*} 
The stabilization term due to artificial matrix can be estimated by Lemma \ref{lemma:stability_stab_term2}, 
\begin{align*}
I_3(\chi) \leq Ch(\|\nabla \theta^{n-1}\|_{L^{2}} + h\|P_hu^{n-1}\|_{2})\|\nabla \chi\|_{L^{2}} \leq C(h^{-1}\|\theta^{n-1}\|_{L^{2}} + h)\|\chi\|_{L^{2}}.
\end{align*}
Similar, the error due to mass lumping,
\begin{align*}
I_4(\chi) \leq Ch^2\|\overline{\partial}\nabla P_hw^{n-1}\|_{L^{2}}\|\nabla \chi\|_{L^{2}} \leq Ch\|\chi\|_{L^{2}}.
\end{align*}
Gathering all these estimates with $\chi = \theta^{n,1} - \theta^{n-1}$ and using Young inequality multiple times, we get
\begin{align*}
\|\theta^{n,1} - \theta^{n-1}\|_{h}^2 \leq Ck^2\,(h^{-2}\|\theta^{n-1}\|_{L^{2}}^2 + h^4 + h^2),
\end{align*}
where since $k = \mathcal{O}(h^2),$ then $k^2\,h^{-2} \leq Ch^2 \leq Ck.$ Hence,
\begin{align*}
\|\theta^{n,1} - \theta^{n-1}\|_{h}^2 \leq Ck(\|\theta^{n-1}\|_{L^{2}}^2 + k\,h^2).
\end{align*}
Using the triangle inequality, we can derive an estimate for $\theta^{n,1},$ i.e.,
\begin{align*}
\|\theta^{n,1}\|_{h}^2 \leq \|\theta^{n,1} - \theta^{n-1}\|_{h}^2 + \|\theta^{n-1}\|_{h}^2 \leq (1+Ck)\|\theta^{n-1}\|_{h}^2 + C\,k^2\,h^2.
\end{align*}
Inserting this estimate into \eqref{estimate_n_prior_final}, we get
\begin{align*}
\|\theta^{n}\|_{h}^2  \leq (1 + Ck)\|\theta^{n-1}\|_{h}^2 + Ck(k^4 + h +\|\theta^{n-1}\|_{h}^2).
\end{align*}
Summing over $n,$ we finally derive the estimate,
\begin{align*}
\|\theta^{n}\|_{h} \leq C(k^2 + h^{1/2}).
\end{align*}
Finally, we obtain the desired estimate by combining the latter estimate together with \eqref{mass_lump_equivalence} and the estimate \eqref{L2_projection_est2_2D}, i.e.,
\begin{align*}
\|U^n - u^n\|_{L^{2}} \leq \|\theta^n\|_{L^{2}} + \|\rho^n\|_{L^{2}} \leq C(k^{2} + h^{1/2}).
\end{align*}

\end{proof}

\section{Numerical experiments}\label{section:numerical_results}
In this section we present several numerical experiments, in order to test the order accuracy of the analyzed fully-discrete scheme, the second order explicit strong stability preserving Runge--Kutta (SSP-RK2) \eqref{fully_disc_1}--\eqref{fully_disc_2}  with respect to the temporal variable.\par
\begin{figure}
\centering
\includegraphics[scale=1]{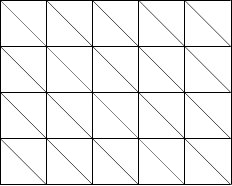}
\caption{A triangulation of a square domain.}\label{fig:triangulation}
\end{figure}

We consider a uniform mesh $\Th$ of the unit square $\Omega = [0,1]^2.$ Each side of $\Omega$ is divided into $M$ intervals of length $h_0=1/M$ for $M\in\mathbb{N}$ and we define the triangulation $\Th$ by dividing each small square by its diagonal, see Fig. \ref{fig:triangulation}. Thus $\Th$ consists of $2M^2$ right-angle triangles with diameter $h = \sqrt{2}h_0.$ 

In order to illustrate the order of convergence of each numerical scheme, we consider the correction factors $\mathfrak{a}_{ij},$ we use Algorithm \ref{algorithm-1} with $q_{i} = \sum_{j\in\Zh^i}d_{ij}(\psi),\,\psi\in\Sh,$ where $d_{ij}$, $i,j=1,\dots,\N$, are the elements of mass matrix $\D_{\psi}.$

\subsection{Maximum principle preservation}
In this subsection, we will study the maximum principle preservation of the stabilized fully-discrete schemes  \eqref{matrix_systemSSP2}. Moreover, we will show that the corresponding initial fully-discrete schemes, i.e., the linear schemes with $\widehat{\mathfrak{a}}_{ij} = 1,\,i,j=1,\ldots,\N$ and $\M_L=\M,$ may not satisfy the discrete maximum principle, \eqref{discrete_max_principle}. We consider the following set of initial conditions for \eqref{conv_law},
\begin{equation}\label{intial_cond_negative}
\begin{aligned}
u^0 & = e^{-100((x-1/2)^2 + (y-1/2)^2)},\;\;(x,y)\in (0,1)^2,\;\;\text{and}\;\;u^0|_{\partial\Omega} = 0,
\end{aligned}
\end{equation}
with final time $T=10k,$ where $k>0,$ will chosen as follows. First, we consider a triangulation of $\Omega,$ as described above with $h_0 = 1/10$ and $k = \frac{1}{10} h^{1.01}.$
In the Table \ref{table:positivity}, we present the values of the coefficient vector of $U^{\NT} \approx u(x,y,T)$ for both standard FEM method and AFC scheme for fixed $y = 0.1$ and $x = [0.0,0.1,\ldots,0.9,1.0].$ Since the initial function is non-negative, then the numerical scheme that preserves the discrete maximum principle should not have negative values. We can see, that these coefficients for this setting can take negative values for the standard FEM scheme while for AFC remains non-negative as expected.

\begin{table}
\begin{center}
\caption{The values of the coefficient vector of $U^{\NT}$ at $x=0,0,1,\ldots,0.9,1,\,y=0.1,$ for the standard FEM and AFC schemes for the solution $u$ with initial conditions \eqref{intial_cond_negative}.}\label{table:positivity}
\begin{tabular}{@{}cccc@{}}
    \toprule
    {$x$}    & {Stand. FEM}         &  {AFC} 	      \\ \hline
    {$0.0$} & {$0.0$}          		& {$0.0$}          \\ \hline
    {$0.1$} & {$1.0309e-02$}        & {$4.6916e-03$}   \\ \hline
    {$0.2$} & {$-2.8695e-03$}  		& {$1.0396e-02$}   \\ \hline
    {$0.3$} & {$-4.1276e-02$}  		& {$1.9539e-02$}   \\ \hline
    {$0.4$} & {$-1.2289e-01$}  		& {$2.9457e-02$}   \\ \hline
    {$0.5$} & {$-2.5104e-01$}  		& {$3.5697e-02$}   \\ \hline
    {$0.6$} & {$-3.4452e-01$}  		& {$3.6743e-02$}   \\ \hline
    {$0.7$} & {$-4.5032e-01$}  		& {$3.2429e-02$}   \\ \hline
    {$0.8$} & {$-2.2527e-01$}  		& {$2.4577e-02$}   \\ \hline
    {$0.9$} & {$-3.9228e-01$}  		& {$2.1874e-02$}   \\ \hline
    {$1.0$} & {$0.0$}          		& {$0.0$}         \\ \bottomrule
\end{tabular}
\end{center}
\end{table}

\subsection{Linear advection equation}
In this section our aim is to numerical illustrate the theoretical findings about convergence rate of the discretization errors of numerical methods \eqref{fully_disc_1}--\eqref{fully_disc_2} in Theorem \ref{theorem:error_estimates_2D}.

\subsubsection{Convergence test for temporal error}
In this subsection, we will study the error convergence for the temporal discretization error of the stabilized fully-discrete scheme \eqref{fully_disc_1}--\eqref{fully_disc_2}. We consider two different initial functions $u^0$ for the linear advection equation, i.e., \eqref{conv_law} with $\bff(u) = \be\,u,$ with $\be=\be(\x)\in W^{1}_{\infty}(\Omega),$ in $\x\in\Omega = [0,1]^2$ and final time $T = 0.1.$ More specifically, the initial functions $u^0$ are given by
\begin{align}
u^0 & = x(1-x)y(1-y),\;\;(x,y)\in [0,1]^2,\label{intial_cond2}\\
u^0 & = \sin(\pi\,x)\sin(\pi\,y),\;\;\;(x,y)\in [0,1]^2.\label{intial_cond1}
\end{align}
The underlying numerical reference solution for each numerical scheme was obtained with $M = 50$ and small time step $k = T/\NT,$ with $\NT = 10^4.$ We compute the approximation for a sequence of $\NT$ as described above with $k = T/\NT$, $\NT = 100,200,400,800,1600,3200$. In Tables \ref{table:AFC_INITIAL_temporal_L2_error_beta_trigonometric_SSP2}--\ref{table:AFC_INITIAL_temporal_L2_error_beta_polynomial_SSP2} we present temporal errors for a fixed mesh step size for the initial functions \eqref{intial_cond1} and \eqref{intial_cond2} with constant flux vector applied to the linear advection equation \eqref{conv_law}, for standard FEM, i.e., the discretization of \eqref{semi_fem_u_2D} via SSP2 for the stabilzed scheme. The order of convergence for both standard FEM-SSP2 and the AFC-SSP2 \eqref{fully_disc_1}--\eqref{fully_disc_2} is two as expected in view of Theorem \ref{theorem:error_estimates_2D}.\par

In addition, we consider initial functions $u^0$ for the linear advection equation, i.e., for \eqref{conv_law} with $\bff(u) = \be\,u,$ where the flux vector is  spatial and time depended, i.e., $\be=\be(x,y,t),$ which is defined as $\be(x,y,t) = ( e^{-t}\sin(\pi\,x), e^{-t}\sin(\pi\,y))^T$ combined with the initial function \eqref{intial_cond1}. The underlying numerical reference solution for each numerical scheme was obtained with $M = 50$ and small time step $k = T/\NT,$ with $\NT = 2000.$ We compute the approximation for a sequence of $\NT$ as described above with $k = T/\NT$, $\NT = 10,20,40,80,160,320.$ 
In Tables \ref{table:AFC_INITIAL_temporal_L2_error_beta_trigonometric_SSP2}--\ref{table:AFC_INITIAL_temporal_L2_error_beta_trigonometric_xyt_SSP2} we present temporal errors for a fixed mesh step size for the initial function \eqref{intial_cond3} with mesh-depended flux vector. Similar to the case where the flux vector is constant, the SSP2 have the optimal order of convergence with respect to the temporal variable.

\begin{table}
\begin{center}
\caption{Temporal $L^{2}-$norm error and convergence order for \eqref{conv_law} with flux function $\bff(u)=\be\,u$ where $\be = (1,3)^T$ and initial function \eqref{intial_cond2}.}\label{table:AFC_INITIAL_temporal_L2_error_beta_trigonometric_SSP2}
\begin{tabular}{@{}cccccc@{}}
    \toprule
    {$\NT$}    & {Stand. FEM}      & {Order}      &  {AFC}	      & {Order}      \\ \hline
    {$100$}    & {$2.3810e-05$}    &              & {$3.7146e-06$}   &            \\ \hline
    {$200$}    & {$5.9110e-06$}    & {$2.0101$}   & {$8.9352e-07$}   & {$2.0556$} \\	\hline
    {$400$}    & {$1.4748e-06$}    & {$2.0029$}   & {$2.1529e-07$}   & {$2.0532$} \\ \hline
    {$800$}    & {$3.6689e-07$}    & {$2.0071$}   & {$5.3229e-08$}   & {$2.0160$} \\ \hline
    {$1600$}   & {$8.9948e-08$}    & {$2.0282$}   & {$1.2939e-08$}   & {$2.0405$} \\\hline
    {$3200$}   & {$2.0715e-08$}    & {$2.1184$}   & {$2.9561e-09$}   & {$2.1299$} \\    \bottomrule
\end{tabular}
\end{center}
\end{table}

\begin{table}
\begin{center}
\caption{Temporal $L^{2}-$norm error and convergence order for \eqref{conv_law} with flux function $\bff(u)=\be\,u$ where $\be = (x^2,2y)^T$ and initial function \eqref{intial_cond1}.}\label{table:AFC_INITIAL_temporal_L2_error_beta_polynomial_SSP2}
\begin{tabular}{@{}cccccc@{}}
    \toprule
    {$\NT$}    & {Stand. FEM}      & {Order}      &  {AFC}	      & {Order}      \\ \hline
    {$100$}    & {$4.9257e-05$}    &              & {$2.3825e-06$}   &            \\ \hline
    {$200$}    & {$1.2295e-05$}    & {$2.0022$}   & {$5.9485e-07$}   & {$2.0019$} \\	\hline
    {$400$}    & {$3.0692e-06$}    & {$2.0022$}   & {$1.5043e-07$}   & {$1.9834$} \\ \hline
    {$800$}    & {$7.6350e-07$}    & {$2.0072$}   & {$3.7275e-08$}   & {$2.0128$} \\ \hline
    {$1600$}   & {$1.8717e-07$}    & {$2.0283$}   & {$9.0855e-09$}   & {$2.0366$} \\\hline
    {$3200$}   & {$4.3104e-08$}    & {$2.1185$}   & {$2.1006e-09$}   & {$2.1128$} \\    \bottomrule
\end{tabular}
\end{center}
\end{table}

\begin{table}
\begin{center}
\caption{Temporal $L^{2}-$norm error and convergence order for \eqref{conv_law} with flux function $\bff(u)=\be\,u$ where $\be(x,y,t) = ( e^{-t}\sin(\pi\,x), e^{-t}\sin(\pi\,y))^T$ and initial function \eqref{intial_cond1}.}\label{table:AFC_INITIAL_temporal_L2_error_beta_trigonometric_xyt_SSP2}
\begin{tabular}{@{}cccccc@{}}
    \toprule
    {$\NT$}  & {Stand. FEM}       & {Order}      &  {AFC}	      & {Order}   \\ \hline
    {$10$}     & {$2.1341e-04$}   &              & {$2.3178e-04$}   &            \\ \hline
    {$20$}     & {$5.3534e-05$}   & {$1.9951$}   & {$5.3745e-05$}   & {$2.1085$} \\	\hline
    {$40$}     & {$1.3404e-05$}   & {$1.9978$}   & {$1.3770e-05$}   & {$1.9646$} \\ \hline
    {$80$}     & {$3.3500e-06$}   & {$2.0004$}   & {$3.4852e-06$}   & {$1.9822$} \\  \hline
    {$160$}    & {$8.3386e-07$}   & {$2.0063$}   & {$8.6949e-07$}   & {$2.0030$} \\  \hline
    {$320$}    & {$2.0448e-07$}   & {$2.0278$}   & {$2.1150e-07$}   & {$2.0395$} \\ \bottomrule
\end{tabular}
\end{center}
\end{table}

\subsubsection{Convergence test for spatial error}

In this subsection, we will study the error convergence for the spatial discretization error of the stabilized fully-discrete schemes \eqref{fully_disc_1}--\eqref{fully_disc_2}. We consider a linear advection equation with constant $\be = (2,4)^T,$ and source term $f(x,y,t),$
\begin{equation}\label{conv_law_system_f}
\begin{cases}
u_t  + \ddiv (\be  u ) = f(x,y,t) , & \text{in }{{\Omega}}\times [0,T],\\
u = 0, & \text{on }\partial {\Omega}\times[0,T],\\
u(\cdot,0)  = u^0, & \text {in }{{\Omega}},
\end{cases}
\end{equation}
with $\Omega = [0,1]^2.$ First, we choose the function $f,$ so as to have as a trigonometric solution, 
\begin{equation}\label{conv_law_system_f_sol}
u(x,y,t) = e^{-t}\sin(\pi\,x)\sin(\pi\,y),\;\;(x,y)\in [0,1]^2,\,\forall\,t\in [0,T].
\end{equation}
Next, we choose the function $f,$ so as to have as a polynomial solution,

\begin{equation}\label{polynomial_system_f_sol}
u(x,y,t) = e^{-t}x(1-x)y(1-y),\;\;(x,y)\in [0,1]^2,\,\forall\,t\in [0,T].
\end{equation}

\begin{table}
\begin{center}
\caption{Spatial $L^{2}-$norm error and convergence order for \eqref{conv_law_system_f} with flux vector $\be = (2,4)^T$ and solution \eqref{conv_law_system_f_sol}.}\label{table:spatial_L2_error_beta_smooth_sol_sinsinSSP2}
\begin{tabular}{@{}ccccccc@{}}
    \toprule
    {$h_0$}  & {Stand. FEM}      & {Order}        &  {AFC}              & {Order}      \\ \hline
    {$1/10$} & {$9.9949e-03$}    &                & {$2.0315e-02$}      &            \\ \hline
    {$1/20$} & {$2.5133e-03$}    & {$1.9916$}     & {$3.9915e-03$}      & {$2.3476$} \\	\hline
    {$1/40$} & {$6.2940e-04$}    & {$1.9976$}     & {$8.0201e-04$}      & {$2.3152$} \\ \hline
    {$1/80$} & {$1.6181e-04$}    & {$1.9597$}     & {$1.7203e-04$}      & {$2.2210$} \\   \bottomrule
\end{tabular}
\end{center}
\end{table}

\begin{table}
\begin{center}
\caption{Spatial $L^{2}-$norm error and convergence order for \eqref{conv_law_system_f} with flux vector $\be = (2,4)^T$ and solution \eqref{polynomial_system_f_sol}.}\label{table:spatial_L2_error_beta_polynomial_sol_sinsinSSP2}
\begin{tabular}{@{}ccccccc@{}}
    \toprule
    {$h_0$}  & {Stand. FEM}      & {Order}      &  {AFC}           & {Order}      \\ \hline
    {$1/10$} & {$7.0302e-04$}    &              & {$1.1442e-03$}   &            \\ \hline
    {$1/20$} & {$1.7652e-04$}	 & {$1.9938$}   & {$2.2970e-04$}   & {$2.3165$} \\	\hline
    {$1/40$} & {$4.4178e-05$}    & {$1.9984$}   & {$4.8254e-05$}   & {$2.2511$} \\ \hline
    {$1/80$} & {$1.1048e-05$}    & {$1.9996$}   & {$1.1045e-05$}   & {$2.1272$} \\   \bottomrule
\end{tabular}
\end{center}
\end{table}

In the Tables \ref{table:spatial_L2_error_beta_smooth_sol_sinsinSSP2}--\ref{table:spatial_L2_error_beta_polynomial_sol_sinsinSSP2}, we present the spatial error for standard FEM, i.e., the discretization of \eqref{semi_fem_u_2D} via SSP2 for the stabilized scheme \eqref{fully_disc_1}--\eqref{fully_disc_2}. We compute the approximation for a sequence of triangulations $\Th$ as described above with $h_0 = 1/M$, $M=10,20,40, 80$. The final time is chosen to be $T = 0.01,$ and we choose $k = \frac{1}{10}h_0.$

\subsection{Inviscid Burger's equation}

In this section our aim is to numerical illustrate the convergence rate of the discretization errors on the numerical method \eqref{fully_disc_1}--\eqref{fully_disc_2} for \eqref{conv_law} with the nonlinear flux function $\bff(u) = \be\,u^2$ with $\be = (1/2,1/2)^T.$ We also perform numerical experiments it the case of mesh and time dependent flux vector $\be,$ that is divergence free, i.e., $\ddiv\be = 0.$

\subsubsection{Convergence test for temporal error}
In this subsection, we will study the error convergence for the temporal discretization error of the stabilized fully-discrete scheme \eqref{fully_disc_1}--\eqref{fully_disc_2}. In addition to the initial functions $u^0$ \eqref{intial_cond1}, we consider also
\begin{align}\label{intial_cond3}
u^0 & = 10e^{-10((x-1/2)^2 + (y-1/2)^2)} + 5,\;\;(x,y)\in (0,1)^2,\;\;\text{and}\;\;u^0|_{\partial\Omega} = 0,
\end{align}
for the inviscid Burger equation in $\x\in\Omega = [0,1]^2$ and final time $T = 0.01.$ For all the numerical experiments presented in this section, the underlying numerical reference solution for each numerical scheme was obtained with $M = 50$ and small time step $k = T/\NT$ with $\NT = 2000.$ We compute the approximation for a sequence of $\NT$ as described above with $k = T/\NT$, $\NT = 10,20,40,80,160,320$. In Tables \ref{table:AFC_INITIAL_temporal_L2_error_beta_trigonometric_SSP2_burger}--\ref{table:AFC_INITIAL_temporal_L2_error_beta_exponential_SSP2_burger_div0_with_time} we present temporal errors for a fixed mesh step size for the initial functions \eqref{intial_cond1} and \eqref{intial_cond3} with constant flux vector applied to the inviscid Burger equation, for standard FEM, i.e., the discretization of \eqref{semi_fem_u_2D} via SSP2 for the stabilzed schemes. The order of convergence for both standard FEM-SSP2 and the AFC-SSP2 \eqref{fully_disc_1}--\eqref{fully_disc_2} is two as expected in view of Theorem \ref{theorem:error_estimates_2D}. 

\begin{table}
\begin{center}
\caption{Temporal $L^{2}-$norm error and convergence order for \eqref{conv_law} with flux function $\bff(u)=\be\,u^2$ with $\be = (1/2,1/2)^T,$ for the initial function \eqref{intial_cond1}.}\label{table:AFC_INITIAL_temporal_L2_error_beta_trigonometric_SSP2_burger}
\begin{tabular}{@{}cccccc@{}}
    \toprule
    {$\NT$}   & {Stand. FEM}      & {Order}      &  {AFC}	      & {Order}      \\ \hline
    {$10$}    & {$1.0427e-04$}    &              & {$1.0816e-04$}   &            \\ \hline
    {$20$}    & {$2.6092e-05$}    & {$1.9986$}   & {$2.7374e-05$}   & {$1.9823$} \\	\hline
    {$40$}    & {$6.5252e-06$}    & {$1.9995$}   & {$6.8675e-06$}   & {$1.9950$} \\ \hline
    {$80$}    & {$1.6299e-06$}    & {$2.0013$}   & {$1.7157e-06$}   & {$2.0010$} \\ \hline
    {$160$}   & {$4.0558e-07$}    & {$2.0067$}   & {$4.2746e-07$}   & {$2.0049$} \\\hline
    {$320$}   & {$9.9445e-08$}    & {$2.0280$}   & {$1.0492e-07$}   & {$2.0265$} \\    \bottomrule
\end{tabular}
\end{center}
\end{table}

\begin{table}
\begin{center}
\caption{Temporal $L^{2}-$norm error and convergence order for \eqref{conv_law} with flux function $\bff(u)=\be\,u^2$ with $\be = (1/2,1/2)^T,$ for the initial function \eqref{intial_cond1}.}\label{table:AFC_INITIAL_temporal_L2_error_beta_polynomial_SSP2_burger}
\begin{tabular}{@{}cccccc@{}}
    \toprule
    {$\NT$}   & {Stand. FEM}      & {Order}      &  {AFC}	      & {Order}      \\ \hline
    {$10$}    & {$1.2388e-09$}    &              & {$1.3402e-09$}   &            \\ \hline
    {$20$}    & {$3.0967e-10$}    & {$2.0001$}   & {$3.3433e-10$}   & {$2.0031$} \\	\hline
    {$40$}    & {$7.7395e-11$}    & {$2.0004$}   & {$8.3472e-11$}   & {$2.0019$} \\ \hline
    {$80$}    & {$1.9326e-11$}    & {$2.0017$}   & {$2.0832e-11$}   & {$2.0025$} \\ \hline
    {$160$}   & {$4.8092e-12$}    & {$2.0067$}   & {$5.1818e-12$}   & {$2.0073$} \\\hline
    {$320$}   & {$1.1800e-12$}    & {$2.0271$}   & {$1.2703e-12$}   & {$2.0283$} \\    \bottomrule
\end{tabular}
\end{center}
\end{table}
 
\begin{table}
\begin{center}
\caption{Temporal $L^{2}-$norm error and convergence order for \eqref{conv_law} with flux function $\bff(u)=\be\,u^2$ with $\be = (1/2,1/2)^T,$ for the initial function \eqref{intial_cond3}.}\label{table:AFC_INITIAL_temporal_L2_error_beta_exponential_SSP2_burger}
\begin{tabular}{@{}cccccc@{}}
    \toprule
    {$\NT$}   & {Stand. FEM}      & {Order}      &  {AFC}	      & {Order}      \\ \hline
    {$10$}    & {$3.4060e-02$}    &              & {$1.3402e-09$}   &            \\ \hline
    {$20$}    & {$8.2307e-03$}    & {$2.0490$}   & {$3.3433e-10$}   & {$2.8642$} \\	\hline
    {$40$}    & {$2.0502e-03$}    & {$2.0052$}   & {$8.3472e-11$}   & {$2.9544$} \\ \hline
    {$80$}    & {$5.1213e-04$}    & {$2.0012$}   & {$2.0832e-11$}   & {$2.9923$} \\ \hline
    {$160$}   & {$1.2748e-04$}    & {$2.0062$}   & {$5.1818e-12$}   & {$2.9994$} \\\hline
    {$320$}   & {$3.1264e-05$}    & {$2.0277$}   & {$1.2703e-12$}   & {$2.9994$} \\    \bottomrule
\end{tabular}
\end{center}
\end{table}

\begin{table}
\begin{center}
\caption{Temporal $L^{2}-$norm error and convergence order for \eqref{conv_law} with flux function $\bff(u)=\be\,u^2,\,\be = (x,-y),$ with initial function \eqref{intial_cond1}.}\label{table:AFC_INITIAL_temporal_L2_error_beta_trigonometric_SSP2_burger_div0}
\begin{tabular}{@{}cccccc@{}}
    \toprule
    {$\NT$}   & {Stand. FEM}      & {Order}      &  {AFC}	      & {Order}      \\ \hline
    {$10$}    & {$1.7011e-04$}    &              & {$1.6178e-04$}   &            \\ \hline
    {$20$}    & {$4.2498e-05$}    & {$2.0010$}   & {$4.0502e-05$}   & {$1.9979$} \\	\hline
    {$40$}    & {$1.0619e-05$}    & {$2.0007$}   & {$1.0109e-05$}   & {$2.0023$} \\ \hline
    {$80$}    & {$2.6515e-06$}    & {$2.0018$}   & {$2.5226e-06$}   & {$2.0027$} \\ \hline
    {$160$}   & {$6.5966e-07$}    & {$2.0070$}   & {$6.2675e-07$}   & {$2.0089$} \\\hline
    {$320$}   & {$1.6173e-07$}    & {$2.0282$}   & {$1.5381e-07$}   & {$2.0267$} \\    \bottomrule
\end{tabular}
\end{center}
\end{table}

\begin{table}
\begin{center}
\caption{Temporal $L^{2}-$norm error and convergence order for \eqref{conv_law} with flux function $\bff(u)=\be\,u^2,\,\be = (e^{-t}\sin(\pi\,y),e^{-t}\sin(\pi\,x)),$ with initial function \eqref{intial_cond3}.}\label{table:AFC_INITIAL_temporal_L2_error_beta_exponential_SSP2_burger_div0_with_time}
\begin{tabular}{@{}cccccc@{}}
    \toprule
    {$\NT$}   & {Stand. FEM}      & {Order}      &  {AFC}	      & {Order}      \\ \hline
    {$20$}    & {$1.0916e-01$}    & 		     & {$2.2936e-02$}   &   \\	\hline
    {$40$}    & {$2.0056e-02$}    & {$2.4443$}   & {$5.0666e-03$}   & {$2.1785$} \\ \hline
    {$80$}    & {$4.8813e-03$}    & {$2.0387$}   & {$1.2702e-03$}   & {$1.9960$} \\ \hline
    {$160$}   & {$1.2097e-03$}    & {$2.0126$}   & {$3.1660e-04$}   & {$2.0043$} \\\hline
    {$320$}   & {$2.9643e-04$}    & {$2.0289$}   & {$7.7503e-05$}   & {$2.0303$} \\    \bottomrule
\end{tabular}
\end{center}
\end{table}

\subsubsection{Convergence test for spatial error}

In this subsection, we will study the error convergence for the spatial discretization error of the stabilized fully-discrete scheme \eqref{fully_disc_1}--\eqref{fully_disc_2}. We consider a inviscid Burger's equation with source term $f(x,y,t),$
\begin{equation}\label{burger_system_f}
\begin{cases}
u_t  + \ddiv (\be\, u^2 ) = f(x,y,t) , & \text{in }{{\Omega}}\times [0,T],\\
u = 0, & \text{on }\partial {\Omega}\times[0,T],\\
u(\cdot,0)  = u^0, & \text {in }{{\Omega}},
\end{cases}
\end{equation}
with $\Omega = [0,1]^2$ and $\be = (1/2,1/2)^T.$ We perform two numerical experiments. In the first, we choose the function $f,$ so as to have the solution $u$ as in \eqref{conv_law_system_f_sol} while in the second experiment, we choose the function $f,$ so as to have the solution $u$ as in \eqref{polynomial_system_f_sol}. In both cases, the final time is chosen $T=0.01.$

\begin{table}
\begin{center}
\caption{Spatial $L^{2}-$norm error and convergence order for \eqref{burger_system_f} with solution \eqref{conv_law_system_f_sol}.}\label{table:spatial_L2_error_burger1_SSP2_trig}
\begin{tabular}{@{}ccccccc@{}}
    \toprule
    {$h_0$}  & {Stand. FEM}      & {Order}          &  {AFC}    		 & {Order}      \\ \hline
    {$1/10$} & {$1.0080e-02$}    &                  &  {$9.9809e-03$}    &            \\ \hline
    {$1/20$} & {$2.5502e-03$}    & {$1.9828$}       &  {$2.5113e-03$}    & {$1.9907$} \\	\hline
    {$1/40$} & {$6.4415e-04$}    & {$1.9852$}   	&  {$6.2900e-04$}    & {$1.9973$} \\ \hline
    {$1/80$} & {$1.6379e-04$}    & {$1.9755$}   	&  {$1.5734e-04$}    & {$1.9992$} \\   \bottomrule
\end{tabular}
\end{center}
\end{table}

\begin{table}
\begin{center}
\caption{Spatial $L^{2}-$norm error and convergence order for \eqref{burger_system_f} with solution \eqref{polynomial_system_f_sol}.}\label{table:spatial_L2_error_burger1_SSP2_polynomial}
\begin{tabular}{@{}ccccccc@{}}
    \toprule
    {$h_0$}  & {Stand. FEM}      & {Order}          &  {AFC}    		& {Order}      \\ \hline
    {$1/10$} & {$7.0302e-04$}    &         	        & {$7.0199e-04$}  	&            \\ \hline
    {$1/20$} & {$1.7652e-04$}    & {$1.9936$}   	& {$1.7638e-04$}   	& {$1.9928$} \\	\hline
    {$1/40$} & {$4.4178e-05$}    & {$1.9984$}   	& {$4.4160e-05$}   	& {$1.9979$} \\ \hline
    {$1/80$} & {$1.1048e-05$}    & {$1.9996$}   	& {$1.1045e-05$}   	& {$1.9993$} \\   \bottomrule
\end{tabular}
\end{center}
\end{table}

In the Tables \ref{table:spatial_L2_error_burger1_SSP2_trig}--\ref{table:spatial_L2_error_burger1_SSP2_polynomial}, we present the spatial error for standard FEM, i.e., the discretization of \eqref{semi_fem_u_2D} via SSP2 for the stabilized scheme \eqref{fully_disc_1}--\eqref{fully_disc_2}. We compute the approximation for a sequence of triangulations $\Th$ as described above with $h_0 = 1/M$, $M=10,20,40, 80$. The final time is chosen to be $T = 0.01,$ and we choose $k = \frac{1}{10}h_0.$

\section{Conclusions}
In this paper, we considered a linear and a nonlinear scalar conservation law on a bounded domain of $\mathbb{R}^2.$ We discretized the spatial using continuous piecewise linear finite elements and we stabilized the semi-discrete scheme via algebraic flux correction method as described in \cite{kuzmin,kuzmin2002,kuzmin2004} and references therein. To compute the correction factors of the algebraic flux correction method, we use a local extremum diminishing flux limiter that is also used in \cite{gabriel4}. The temporal variable were discretized by the second order strong stability preserving Runge--Kutta method. Under assumptions for the triangulation used for the space discretization and the time step, we derived error estimates in $L^{2}-$norm in space and $\ell^\infty$ in time for the fully-discrete scheme. Numerical experiments in two dimensions were presented for both the standard FEM and stabilized schemes that validates the theoretical results for temporal order of convergence. The spatial order of convergence in the numerical examples is the optimal that it can be achieved when continuous piecewise linear finite elements is used, i.e., the $L^{2}$ spatial error is proportional to $h^2$, while in theoretical part it can be proved that the spatial error in $L^{2}$ is proportional to $h^{1/2}.$

\section*{Acknowledgments}
The author would like to thank Gero Schnücke for contributing several valuable remarks.

\bigskip
\bibliographystyle{acm}
\bibliography{ref}

\begin{thebibliography}{10}

\bibitem{gabriel1}
{\sc Barrenechea, G.~R., John, V., and Knobloch, P.}
\newblock Analysis of algebraic flux correction schemes.
\newblock {\em SIAM J. Numer. Anal 54\/} (2016), 2427--2451.

\bibitem{gabriel4}
{\sc Barrenechea, G.~R., John, V., and Knobloch, P.}
\newblock An algebraic flux correction scheme satisfying the discrete maximum
  principle and linearity preservation on general meshes.
\newblock {\em Math. Models Methods Appl. Sci. 27\/} (2017), 525–548.

\bibitem{barrenechea2024}
{\sc Barrenechea, G.~R., John, V., and Knobloch, P.}
\newblock Finite element methods respecting the discrete maximum principle for
  convection-diffusion equations.
\newblock {\em SIAM Rev. 66}, 1 (2024), 3--88.

\bibitem{gabriel3}
{\sc Barrenechea, G.~R., John, V., Knobloch, P., and Rankin, R.}
\newblock A unified analysis of algebraic flux correction schemes for
  convection–diffusion equations.
\newblock {\em SeMA 75\/} (2018), 655–685.

\bibitem{boman2006}
{\sc Boman, M.}
\newblock Estimates for the ${{L^2}}$-projection onto continuous finite element
  spaces in a weighted ${{L^p}}$-norm.
\newblock {\em BIT Numer. Math. 46}, 2 (2006), 249--260.

\bibitem{brenner}
{\sc Brenner, S.~C., and Scott, L.~R.}
\newblock {\em The Mathematical Theory of Finite Element Methods}, second~ed.
\newblock Springer, New York, 2008.

\bibitem{burman2015}
{\sc Burman, E.}
\newblock A monotonicity preserving, nonlinear, finite element upwind method
  for the transport equation.
\newblock {\em Appl. Math. Lett. 49\/} (2015), 141--146.

\bibitem{burman2010}
{\sc Burman, E., Ern, A., and Fernández, M.~A.}
\newblock Explicit {{R}}unge--{{K}}utta schemes and finite elements with
  symmetric stabilization for first-order linear {{PDE}} systems.
\newblock {\em SIAM J. Numer. Anal. 48}, 6 (2010), 2019--2042.

\bibitem{chatzipa2}
{\sc Chatzipantelidis, P., Lazarov, R., and Thomee, V.}
\newblock Some error estimates for the lumped mass finite element method for a
  parabolic problem.
\newblock {\em Math. Comput. 81\/} (2012), 1--20.

\bibitem{chatzipa3}
{\sc Chatzipantelidis, P., and Pervolianakis, C.}
\newblock Error analysis of a backward {{E}}uler positive preserving stabilized
  scheme for a {{C}}hemotaxis system.
\newblock {\em arXiv preprint arXiv:2210.04709\/} (2022).

\bibitem{cockburn1996}
{\sc Cockburn, B., and Gremaud, P.-A.}
\newblock Error estimates for finite element methods for scalar conservation
  laws.
\newblock {\em SIAM J. Numer. Anal. 33}, 2 (1996), 522--554.

\bibitem{cockburn1989}
{\sc Cockburn, B., and Shu, C.-W.}
\newblock {{TVB}} {{R}}unge--{{K}}utta local projection discontinuous
  {{G}}alerkin finite element method for conservation laws. ii. general
  framework.
\newblock {\em Math. Comput. 52}, 186 (1989), 411--435.

\bibitem{crouzeix1987}
{\sc Crouzeix, M., and Thomée, V.}
\newblock The stability in ${{L^p}}$ and ${{W^{1}_p}}$ of the
  ${{L^2}}$-projection onto finite element function spaces.
\newblock {\em Math. Comput. 48}, 178 (1987), 521--532.

\bibitem{dafermos2005}
{\sc Dafermos, C.~M.}
\newblock {\em Hyperbolic Conservation Laws in Continuum Physics}.
\newblock Springer Verlag, 2nd edition, 2005.

\bibitem{gottlieb2001}
{\sc Gottlieb, S., Shu, C.-W., and Tadmor, E.}
\newblock Strong stability-preserving high-order time discretization methods.
\newblock {\em SIAM Rev. 43}, 1 (2001), 89--112.

\bibitem{guermond2014}
{\sc Guermond, J.-L., and Nazarov, M.}
\newblock A maximum-principle preserving ${{\mathcal{C}^0}}$ finite element
  method for scalar conservation equations.
\newblock {\em Comput. Methods Appl. Mech. Eng. 272\/} (2014), 198--213.

\bibitem{guermond2016}
{\sc Guermond, J.-L., and Popov, B.}
\newblock Invariant domains and first-order continuous finite element
  approximation for hyperbolic systems.
\newblock {\em SIAM J. Numer. Anal. 54}, 4 (2016), 2466--2489.

\bibitem{guermond2017}
{\sc Guermond, J.-L., and Popov, B.}
\newblock Invariant domains and second-order continuous finite element
  approximation for scalar conservation equations.
\newblock {\em SIAM J. Numer. Anal. 55}, 6 (2017), 3120--3146.

\bibitem{hajduk2023}
{\sc Hajduk, H., and Rupp, A.}
\newblock Analysis of algebraic flux correction schemes for semi-discrete
  advection problems.
\newblock {\em BIT Numer. Math. 63}, 1 (2023), 8.

\bibitem{douglas1975}
{\sc Jr., J.~D., Dupont, T., and Wahlbin, L.}
\newblock The stability in ${{L^q}}$ of the ${{L^2}}$-projection into finite
  element function spaces.
\newblock {\em Numer. Math. 23\/} (1975), 193--197.

\bibitem{kuzmin}
{\sc Kuzmin, D.}
\newblock {\em A Guide to Numerical Methods for Transport Equations}.
\newblock University Erlangen-Nuremberg, Nuremberg, 2010.

\bibitem{kuzmin2002}
{\sc Kuzmin, D., and Turek, S.}
\newblock Flux correction tools for finite elements.
\newblock {\em J. Comput. Phys. 175}, 2 (2002), 525--558.

\bibitem{kuzmin2004}
{\sc Kuzmin, D., and Turek, S.}
\newblock High-resolution {{FEM-TVD}} schemes based on a fully multidimensional
  flux limiter.
\newblock {\em J. Comput. Phys. 198\/} (2004), 131–158.

\bibitem{kucera2016}
{\sc Kučera, V.}
\newblock Finite element error estimates for nonlinear convective problems.
\newblock {\em J. Numer. Math. 24}, 3 (2016), 143--165.

\bibitem{kucera2019}
{\sc Kučera, V., and Shu, C.-W.}
\newblock On the time growth of the error of the {{DG}} method for advective
  problems.
\newblock {\em IMA J. Numer. Anal. 39}, 2 (2019), 687--712.

\bibitem{leveque2012}
{\sc LeVeque, R.~J.}
\newblock {\em Numerical Methods for Conservation Laws}, 2~ed.
\newblock Lectures in Mathematics. ETH Zürich. Birkhäuser Basel, 2012.

\bibitem{lohmann}
{\sc Lohmann, C.}
\newblock Algebraic flux correction schemes preserving the eigenvalue range of
  symmetric tensor fields.
\newblock {\em ESAIM: M2AN 53\/} (2019), 833–867.

\bibitem{ciarlet}
{\sc Philippe, C.}
\newblock {\em The Finite Element Method for Elliptic Problems}.
\newblock Society for Industrial and Applied Mathematics, North-Holland, 2002.

\bibitem{zhang2004}
{\sc Zhang, Q., and Shu, C.}
\newblock Error estimates to smooth solutions of {{R}}unge--{{K}}utta
  discontinuous {{G}}alerkin methods for scalar conservation laws.
\newblock {\em SIAM J. Numer. Anal. 42\/} (2004), 641--666.

\bibitem{zhang2010}
{\sc Zhang, Q., and Shu, C.}
\newblock Stability analysis and a priori error estimates of the third order
  explicit {{R}}unge--{{K}}utta {{D}}iscontinuous {{G}}alerkin {{M}}ethod for
  scalar conservation laws.
\newblock {\em SIAM J. Numer. Anal. 48}, 3 (2010), 1038--1063.

\bibitem{zhang2010_2}
{\sc Zhang, X., and Shu, C.-W.}
\newblock On maximum-principle-satisfying high order schemes for scalar
  conservation laws.
\newblock {\em J. Comput. Phys. 229}, 9 (2010), 3091--3120.

\end{thebibliography}

\end{document}